\magnification=\magstep1 \overfullrule=0pt \def\ga{\gamma}  \def\la{\lambda}
\def\r{{\overline{r}}} \def\k{{\overline{k}}} \def\R{{\cal
R}}\def\NZ{{\cal NZ}}
\def\qed{\quad \vrule height 1.9ex width 1.3ex depth -.1ex}
\def\i{{\rm i}} \def\G{\Gamma}  \def\si{\sigma}



\def\Z{\ifmmode{Z\hskip -4.8pt Z}
     \else{\hbox{$Z\hskip -4.8pt Z$}}\fi}
\def\Q{{\ifmmode{Q\hskip-5.0pt\vrule height6.0pt depth 0pt\hskip6pt}
      \else{\hbox{$Q\hskip-5.0pt\vrule height6.0pt depth
0pt\hskip5pt$}}\fi}}
\def\C{{\ifmmode{C\hskip-5.0pt\vrule height6.0pt depth 0pt\hskip6pt}
      \else{\hbox{$C\hskip-5.0pt\vrule height6.0pt depth
0pt\hskip5pt$}}\fi}}

\font\small=cmr8
\font\smal=cmr7
\font\smcap=cmcsc10
\hoffset=1truecm\hsize=15truecm   \vsize=22truecm
\baselineskip=15pt
{\nopagenumbers
\rightline{DAMTP-97-101} 
\rightline{September, 1997}\vskip1cm

\centerline{{\bf ON FUSION ALGEBRAS AND MODULAR MATRICES}\footnote{$^*$}{\small 
This
research  was supported in part by NSERC.}}

\vskip1cm\centerline{T. Gannon\  and M.A. Walton\
\footnote{$^\dagger$}{\small On leave from the Physics Dept, Univ. Lethbridge, Alberta,
Canada}} \vskip.5cm

\centerline{{\it Mathematics Department,
York University}} \centerline{{\it North York, Ontario,
Canada\ \ M3J 1P3}}\centerline{\rm tgannon@mathstat.yorku.ca}

\vskip.25cm
\centerline{\it Department of Applied Mathematics and Theoretical Physics}
\centerline{\it 
University of Cambridge} \centerline{{\it Silver Street,
Cambridge\ \ CB3 9EW, U.K.}}\centerline{\rm m.a.walton@damtp.cam.ac.uk} 

\vskip 1truecm  
\noindent{{\bf Abstract:} We consider the fusion algebras arising in e.g.\
 Wess-Zumino-Witten conformal field theories, affine Kac-Moody algebras
at positive integer level, and quantum groups at roots of unity. Using 
properties of the modular
matrix $S$, we find small sets of primary fields (equivalently, sets
of highest weights) which can be identified with
the variables of a polynomial realization of the $A_r$ fusion algebra at level
$k$.
We prove that for many choices of rank $r$ and level $k$, the number of these
variables is the minimum possible, and we conjecture that it is in fact
minimal for most $r$ and $k$. We also find new, systematic sources of
zeros in the modular matrix $S$. In
addition, we obtain a formula relating the entries of $S$ at
fixed points, to entries of $S$ at smaller ranks and levels. Finally,
we identify the number fields generated over the rationals by the entries of 
$S$, and by the fusion (Verlinde) eigenvalues.} \vskip1cm

\vfill\eject}
\pageno=1

\noindent{\bf 1. Introduction} \medskip

Fix an affine non-twisted algebra $g=X_r^{(1)}$, and level $k$. Put $\k:=
k+h^\vee$, where $h^\vee$ is the dual Coxeter number of $g$. Let $w^0,\ldots,
w^r$ denote its fundamental weights, and put $\rho:=\sum_{i=0}^rw^i$. Let
$P_+^k(g)$ be the set of all level $k$ integrable highest weights of $g$.
For example,
$$P_+^k(A_r^{(1)})=\{\sum_{i=0}^r \la_iw^i\,|\,\la_i\in\Z_{\ge 0},\
\sum_{i=0}^r\la_i=k\}\ .$$
Write $ch_\la$ for the corresponding character. Sometimes it is convenient
to write $(\la_0,\la_1,\ldots,\la_r)$ for $\sum_i\la_i w^i$.
When the level of a weight
is known, we will often drop the $w^0$ component. For example,
the element $kw^0$ of $P_+^k(g)$ will be denoted by 0.
The corresponding quantities for the underlying finite-dimensional Lie
algebra $\bar{g}$ will always be denoted with a bar.

Under the familiar action of $SL_2(\Z)$ on the Cartan subalgebras of $g$,
we find that the span of the level $k$ characters $ch_\la$ is stable. In
particular,  define a matrix $S$ by:
$$ch_\la\left({-1\over \tau},{z\over \tau},u- {(z|z)\over \tau}\right)=\sum_{\mu\in
P_+^k(g)}   S_{\la,\mu}\,ch_\mu(\tau,z,u)\ .$$
$S$ has several interesting properties. Most importantly:

\medskip
\noindent{{\smcap Lemma 1}} (Kac-Peterson [16]):\quad {\it Let $\overline{ch}_{\bar{
\nu}}$
denote the Weyl character of $\bar{g}$ with highest weight $\bar{\nu}$. Then
for any $\la,\mu\in P_+^k(g)$, we have both $S_{0,\mu}\ne 0$ and}
$${S_{\la,\mu}\over S_{0,\mu}}=\overline{ch}_{\overline{\la}}\left(-2\pi \i
{\overline{\mu+\rho}\over \k}\right)=:\chi_\la(\mu)\ .\eqno(1.1a)$$
\medskip

By Lemma 1, a useful expression for $\chi_\la(\mu)$ is
$$\chi_\la(\mu)=
\sum_{\beta\in\overline{\Omega}(\overline{\la})}\sum_{\ga\in\overline{W}
\beta}m_{\overline{\la}}(\beta)\,\exp[-2\pi \i\,{ \ga\cdot(\overline{\mu+\rho})
\over\k}]\ ,\eqno(1.1b)$$
where $\overline{W}$ is the (finite) Weyl group, where $\overline{\Omega}(
\overline{\la})$ is the set of dominant weights of the representation of
$\bar{
g}$ with highest weight $\overline{\la}$, and where $m_{\overline{\la}}(\beta)$
is the weight multiplicity.

\vfill\eject A classical result is:

\medskip\noindent{{\smcap Lemma 2}} (Cartan [3]): \quad {\it For each
$\bar{\nu}$, we can write $\overline{ch}_{\bar{\nu}}=P_{\bar{\nu}}
(\overline{ch}_{\overline{w^1}},\ldots,\overline{ch}_{\overline{w^r}})$ for
some polynomial  $P_{\bar{\nu}}(x_1,\ldots,x_r)$.}\medskip

Therefore,
$$\chi_\la(\mu)=P_{\bar{\la}}(\chi_{w^1}(\mu),\ldots,\chi_{w^r}(\mu))\ ,
\eqno(1.2)$$
for all $\mu\in P_+^k(g)$.

Define the {\it fusion matrices} $N_\la$ by Verlinde's formula [21]:
$$(N_\la)_\mu^\nu:=N_{\la,\mu}^\nu=\sum_{\ga\in P_+^k(g)}S_{\la,\ga}
{S_{\mu,\ga}\over S_{0,\ga}}S_{\nu,\ga}^*\ .\eqno(1.3)$$
Equation (1.3) tells us that the $N_\la$ are simultaneously diagonalized by $S$,
and have eigenvalues $\chi_\la(\mu)$.
The {\it fusion algebra} (or {\it Verlinde algebra}) of $g$ at level $k$ is
defined to be the $\C$-span of
$\{N_\la\ :\ \la\in P_+^k(g)\}$. It is associative and commutative, with
 unit $N_0=I$ and  integer structure constants $N_{\la,\mu}^\nu$:
$$N_\la\,N_\mu=\sum_{\nu\in P_+^k(g)}N_{\la,\mu}^\nu\,N_\nu\ .$$
In fact it is isomorphic as an algebra to $\C^{\|P_+^k(g)\|}$, defined with
componentwise addition and multiplication, and so a critical ingredient here
in our definition is the choice of preferred basis $\{N_{\la}:\la\in P_+^k(g)\}$.
Fusion algebras (or the corresponding fusion ring) appear in many
different contexts, e.g.\ in rational conformal
field theory (RCFT) [21]. The RCFTs with fusion
algebras of the type discussed here, i.e.\ those associated with some $g$, are
known as Wess-Zumino-Witten models. Fusion algebras also appear in
the study of quantum groups [19] and Hecke algebras [14] at roots of
unity, Chevalley groups at nonzero characteristic [12], and quantum
cohomology [22].

Call a set $\Gamma=\{\ga^1,\ldots,\ga^n\}\subset P_+^k(g)$ a
{\it fusion-generator}
if any $N_\la$ can be written as a polynomial\footnote{$^1$}{{\smal By
Lagrange interpolation, `polynomial' here is equivalent to `function'.}}
in $N_{\ga^1},\ldots,N_{\ga^n}$ -- in other words, if for each $\la\in P_+^k(
g)$ there is a polynomial $P_\la(x_1,\ldots,x_n)$ such that
$$\chi_\la(\mu)=P_\la(\chi_{\ga^1}(\mu),\ldots,\chi_{\ga^n}(\mu))\qquad
\forall \mu\in P_+^k(g).\eqno(1.4a)$$
Equivalently, $\Gamma$ is a fusion-generator\footnote{$^2$}{{\smal Our
definition should not be confused with the `bootstrapped' version of
a fusion-generator used in [10].}} iff for any $\la,\mu\in
P_+^k(g)$, the only way we can have
$$\chi_{\ga^\ell}(\la)=\chi_{\ga^\ell}(\mu)\quad{\rm for\ all}\ \ell=1,
\ldots,n,\eqno(1.4b)$$
is when $\la=\mu$. 

The equivalence of  the statements of (1.4a) and (1.4b) can be
seen as follows. First, if (1.4b) holds, then (1.4a) implies
$\chi_\phi(\lambda)=\chi_\phi(\mu)$ for all $\phi\in
P_+^k(g)$. Multiplying this result by $S^*_{\nu,\phi}$ and summing
over $\phi\in P_+^k(g)$ gives $\lambda=\mu$, by the unitarity of the
matrix $S$. 

In the other direction, we need to construct a polynomial $P_\lambda$ 
in $n=\|\Gamma\|$ variables, taking values $\chi_\lambda(\mu)$ at 
$m=\|P_+^k(g)\|$ distinct points.
Let $\vec x:=(x_1,\ldots,x_n)$ denote a point in 
$\C^n$, and let $\vec x_a$, $a=1,\ldots,m$ be the points at which the required
polynomial must take the values $y_a$. Here $x_{a,j}=
\chi_{\gamma^j}(\mu^a)$ and $y_a=\chi_\la(\mu^a)$, where $a$ labels the
different weights of
$P_+^k(g)$. A polynomial of minimal degree satisfying the requirements
can be constructed by the Lagrange interpolation formula:
$$P(\vec x) =
\sum_{a=1}^m y_a\prod^m_{b=1,b\not=a} {\vec r\cdot(\vec x-\vec
x_b)\over\vec r\cdot(\vec x_a-\vec x_b)}\ .$$
Here $\vec r$ can be any (constant) vector such that
$\vec r\cdot(\vec x_a-\vec x_b)$ vanishes iff $a=b$.  

By the {\it fusion-rank} $\R_k(g)$, we mean the minimum
possible cardinality $n=\|\Gamma\|$ of a fusion-generator $\Gamma$. Such
a $\Gamma$ is called a {\it fusion-basis}.

\medskip\noindent{{\bf Question 1:}}\quad For a given $g$ and $k$,
what is the fusion-rank $\R_k(g)$, and what is a fusion-basis?\medskip

This problem was studied by Di Francesco and Zuber [6]. For the
applications it should suffice to get a reasonable upper bound for the
fusion-rank, and to find a $\Gamma$ which realizes that bound.
Incidently, it was proven in [1] that there will be a fusion potential [13]
corresponding to any fusion-generator $\Gamma$.

Question 1 seems a natural one from the fusion algebra
perspective, and is especially interesting considering that the fusion-rank
often turns out to be surprisingly low. This analysis should have consequences
for the work of Moody, Patera, Pianzola, $\ldots$ on elements of finite order in a
finite-dimensional Lie group (see e.g.\ [18,20] and references therein). It has
direct relevance for the classification of conformal field
theories (more precisely, their 1-loop partition functions; see e.g.\
[9,11,10]). Our results may lead to a
new presentation of the fusion algebras, along the lines of the Schubert
calculus of [13,15]. As another example, we
mention that our problem may be related to finding bases for the quantum
cohomology of Grassmannians [22]. 

Incidentally, these fusion algebras all have a rank of one, in a sense:
precisely, the Krull dimension of a fusion algebra 
will be one. It is not difficult to find an element $N$ of the fusion algebra
in which every fusion matrix $N_{\la}$ will be a polynomial. These $N$ however
will in general be nontrivial linear combinations of our basis vectors (1.3).
For the applications we are interested in, this observation is not helpful.
There is a natural basis for the fusion algebra, namely $P_+^k(g)$,
and an important condition is that fusion-generators
are required to be subsets of that basis.

We will address Question 1 for $g=A_r^{(1)}$ in Section 3. Our
best lower bound for $\R_{r,k}:=\R_k(A^{(1)}_r)$ is given in Thm.\ 1(2);  our
best upper
bound and smallest fusion-generator is given in Thm.\ 3. Cor.\ 1 tells
us precisely
when $\{w^1\}$ is a fusion-generator. Cor.\ 2 answers Question
1 when $r$ or $k$ is small, and Conjecture 1 gives our guess for a 
general statement.

Another
question related to this one, which we will consider in Section 4, is:

\medskip\noindent{{\bf Question 2:}}\quad For $g=A_r^{(1)}$, when is $N_{w^1}$
invertible?\medskip

The first fundamental weight $w^1$ is especially interesting, since (1.1b)
and its fusion numbers $N_{w^1,\mu}^\nu$ are so simple. Incidentally, 
$N_\la$ is invertible iff $N_{\la^\sigma}$ is, for any Galois element $\sigma$
(see (2.6) below) -- this holds in fact for any RCFT [5]. However, the
inverse of a fusion matrix will only
itself be a fusion matrix in the trivial cases: $(N_\la)^{-1}=N_\mu$ iff both
$\la=J^a0$ and $\mu=J^{-a}0$ for some $a\in\Z$, where $J$ is given in (2.1b) --
again the analogue holds for any RCFT. (The proof of this uses the fact
that the inverse of a non-negative integer matrix can itself be integral
and non-negative, only if it is a permutation matrix.)

Our best condition for $N_{w^1}$ being invertible is given in Thm.\ 6(3), while
our best conditions for noninvertibility are Thms.\ 6(4),(5). Together,
these answer Question 2 for most $r,k$. Conjecture 2
gives our guess for the general answer.

A final question, which we solve in Section 6, was asked in [4]. It is interesting
because of the Galois action (2.6) on the matrix $S$ and on the fusion coefficients.

\medskip\noindent{{\bf Question 3:}}\quad For $A_r^{(1)}$, what are the number
 fields $K_{r,k}$ and $L_{r,k}$ generated over the rationals by the entries
 $S_{\la,\mu}$, and by the fusion (Verlinde) eigenvalues $\chi_\la(\mu)$, respectively?

\bigskip\noindent{{\bf 2. The $A_{r,k}$ Modular Matrix $S$}}\medskip

For now, let us restrict attention to $A_{r,k}$ (i.e.\ $A_r^{(1)}$ at level
$k$). Write $\r:=r+1$, $P_+^{r,k}:=
P_+^k(A_r^{(1)})$ and $\R_{r,k}:=\R_k(A_r^{(1)})$. The symmetry group of its
Coxeter-Dynkin diagram is the dihedral group on $\r$ elements, generated by
an order 2  {\it conjugation} $C$ and an order $\r$ {\it simple current} $J$:
$$\eqalignno{C\la=&\,\la_0w^0+\sum_{i=1}^r\la_{r+1-i}w^i\ ,&(2.1a)\cr
J\la=&\,\la_rw^0+\sum_{i=1}^r\la_{i-1}w^i\ .&(2.1b)\cr}$$
These act on the $\chi_\la(\mu)$ by
$$\eqalignno{\chi_{C\la}(\mu)=&\,\chi_\la(C\mu)=\chi_\la(\mu)^*\ ,&(2.2a)\cr
\chi_{J^a\la}(J^b\mu)=&\,\exp[2\pi \i(b\,t(\la)+a\,t(\mu)+kab)/\r]\,
\chi_{\la}(\mu)\ , &(2.2b)\cr}$$
where
$$t(\la):=\sum_{j=1}^rj\la_j\eqno(2.2c)$$
is called the $\r$-ality. A useful relation is
$$t(J^a\la)\equiv ak+t(\la)\quad ({\rm mod}\ \r)\ .\eqno(2.2d)$$ 

Another `symmetry' of $\chi_\la(\mu)$, when $k\ne 1$, is {\it rank-level
duality} [2]:
$$\chi_\la(\mu)=\exp[2\pi \i\, t(\la)\,t(\mu)/\r k]\,
\widetilde{\chi}_{\tau\la}(\tau\mu)^*\ ,\eqno(2.3a)$$
where $\tau\la$ denotes the weight in $P_+^{k-1,r+1}$ corresponding to
the transpose (sometimes called `conjugate') of the Young diagram of $\la$,
after deleting any columns of
length $k$ in the transposed diagram (reminder: the $i$th row of the Young
diagram of $\la$ has $\sum_{j=i}^r\la_j$ boxes). This 
deletion is a consequence of (2.4f) below. We will usually denote the
quantities of $A_{k-1,r+1}$ with tildes.
For example, $\tau w^\ell=\ell\widetilde{w}^1$. $\tau$ defines a bijection
between
the $J$-orbits in $P_+^{r,k}$ and the $\tilde{J}$-orbits in $P_+^{k-1,r+1}$.
Note that
$$\tilde{t}(\tau\la)\in t(\la)-k\Z_{\ge 0}\ ,\eqno(2.3b)$$
since $t(\la)$ is the number of boxes in the Young diagram of $\la$. 

The Weyl group of $A_r$ is the symmetric group $S_{\r}$. This gives us an
essential property of $S$: its relation to the symmetric polynomials. In
particular, we can see from (1.1b) that
$$\chi_\la(\mu)=\exp[2\pi \i\,t(\la)\,t(\mu+\rho)/\r\k]\,S_\la(x_1,\ldots,
x_{\r})\ ,\eqno(2.4a)$$
where $x_\ell:=\exp[-2\pi \i\,\mu(\ell)/\k]$ for $\mu(\ell):=\sum_{j=\ell}^r
(\mu_j+1)$. $S_\la$
is a polynomial over $\Z$ -- the {\it Schur polynomial} of shape $(\sum_{i=1}^r
\la_i,\sum_{i=2}^r\la_i,\ldots,\la_r)$
[8] -- symmetric in the $x_i$, and homogeneous of degree $t(\la)$. It is
often convenient to write $S_\la$ as a polynomial
$$Q_\la(y_1,\ldots,y_{\r k})=\sum_{m=(m_1,\ldots,m_{\r k})} c_{m}\prod_\ell
y_\ell^{m_\ell}\ ,\eqno(2.4b)$$
evaluated at the `power sums' of our $x_i$:
$$y_\ell=\sum_{i=1}^{\r}
x_i^\ell=P_\ell(x_1,\ldots,x_{\r})\ .\eqno(2.4c)$$
The coefficients $c_m$ of $Q_\la$ can be expressed in terms of the characters
of the symmetric group $S_{\r}$ (this is essentially Frobenius-Schur duality),
and each nonzero  $c_{m}$  will have
$\sum jm_j=t(\la)$ [8]. We will also write $S_\la[\mu]$ and $P_\ell[\mu]$, when
convenient. Note that
$$P_\ell[J^m\mu]=\exp[2\pi \i\,\ell\,\mu(\r-m)/\k]\,P_\ell[\mu]\ .\eqno(2.4d)$$
A valuable special case of (2.4a) is
$$\chi_{w^\ell}(\mu)=\exp[2\pi \i\,\ell\,t(\mu+\rho)/\r \k]\sum_{1\le i_1<
\cdots<i_\ell\le \r}x_{i_1}\cdots x_{i_\ell}\ .\eqno(2.4e)$$

Symmetric polynomials have an important variable-specialisation property
which permits the number of variables to be increased (with the extra variables
set to 0), and yet {\it all algebraic relations\footnote{$^3$}{{\smal such as
(2.4), but not e.g.\ (2.3a), (2.6) or (2.8). More precisely, specialisation
defines a homomorphism between the polynomial rings, taking Schur polynomials
to Schur polynomials, power sums to power sums, etc.}}
among the symmetric polynomials
will be preserved}. This permits us to define $\chi_{\la}$ when $\la$ has more
than $\r$ components, using (2.4a) with variables $x'_1=x_1,\ldots,
x'_\r=x_\r$, and $x'_{\r+1}=\cdots=0$, and we find
$$\chi_{(\la_0,\la_1,\ldots,\la_{r},\ldots)}(\mu)=\left\{\matrix{0&{\rm if}\
\la_\ell>0\ {\rm for\ some}\ \ell>\r\cr \chi_{(\la_0,\la_1,\ldots,\la_r)}(\mu)
&{\rm otherwise}\cr}\right. ,\eqno(2.4f)$$
valid for any $\mu\in P_+^{r,k}$. This can be directly understood using
for example the construction of Schur polynomials from Young Tableaux.
A special case of (2.4f) is $\chi_{w^{\r}}=1$
and $\chi_{w^\ell}=0$ for $\ell>\r$. We will use (2.4f) in several places --
see e.g.\ the proof of Thm.\ 3. 

Call $\la\in P_+^{r,k}$ a $J^d$-{\it fixed point} if $d$ is the
smallest positive integer satisfying $J^d\la=\la$ -- in other words if the
$\la_i$ have period $d$. We will say $\la$ is a
{\it fixed point} if it is a $J^d$-fixed point for some $d<\r$. 
Note that if $\varphi$ is a fixed point of $J^d$, 
we can speak of a `truncated weight'
$(\varphi_0,\varphi_1,\ldots,\varphi_{d-1})=:\varphi^\prime$; by (2.5a) below
it will lie in $P_+^{d,kd/\r}$. We have 
$$\eqalignno{{dk\over \r}=&\,\sum_{i=0}^{d-1}\varphi_i=\,\sum_{i=0}^{d-1}
\varphi^\prime_i&(2.5a)\cr
t(\varphi)=&\,{\r\over d}\sum_{j=1}^{d-1}j\varphi_j+k\,{\r-d\over 2}
=\,{\r\over d}t^\prime(\varphi^\prime)+k\,{\r-d\over 2}\ ,&(2.5b)\cr}$$
where $t'$ denotes $d$-ality.
There exist $J^d$-fixed points in $P_+^{r,k}$ iff $d$ divides $\r$ and $\r/d$
divides $k$. In other words, the smallest fixed-point period is $\r/{\rm gcd}
\{\r,k\}$, and all other possible periods are multiples of this number.
Also, if $\varphi$ is a $J^{\r/d}$-fixed point, its
rank-level dual $\tau\varphi$ is a $\tilde J^{k/d}$-fixed point.   

By (2.2b), if $\mu$ is a
$J^d$-fixed point, then $\chi_{\la}(\mu)=0$
whenever $t(\la)\not\equiv 0$ (mod $\r/d$). The same comment holds for
$\mu$ if
instead $\la$ is a $J^d$-fixed point. This is certainly not the
only source of zeros in the  matrix $S$ however, as we shall see, but it
is an important one. In fact, there are {\it many} more zeros at fixed points
than this simple $\r$-ality test suggests. For example, of all weights $\la$
with $t(\la)=\r/d$, the entry $S_{\la,\varphi}$ will equal zero for {\it every}
$J^d$-fixed point $\varphi$, unless $\la$ is a hook $({\r\over d}-a)w^1+w^a$.
We will describe below the set ${\cal NZ}(d)$ of all weights $\la$
which can have nonzero entries at $J^d$-fixed points.

Moreover, many different weights $\la\ne\mu$ -- even in the set ${\cal NZ}
(d)$ -- will have the same value $S_{\la,\varphi}=S_{\mu,\varphi}$ at
all $J^d$-fixed points 
$\varphi$. For example, for the hooks $\la$ with $t(\la)={\r\over d}$,
we will have $\chi_{\la}(\varphi)=\pm \chi_{w^{\r/d}}(\varphi)$ for all
$\varphi$, where the sign is independent of $\varphi$. More generally, note that
the right side of (2.8c) is independent of $a''$, except for the unimportant
sign.

Hence fixed point considerations are very important for both Questions 1 and 2,
and play a large role in this paper.

An unexpected symmetry of the matrix $S$ is the Galois action discussed in
[5]. For any $\si\in{\rm Gal}(K_{r,k}/\Q)$, there exists a
permutation $\mu\mapsto \si \mu$ of 
$P_+^{r,k}$ such that   
$$\eqalignno{\si S_{\la,\mu}=&\,\epsilon_\si(\mu)\,S_{\la,\si\mu}\ ,&(2.6a)\cr
\si\chi_\la(\mu)=&\,\chi_\la(\si \mu)\ ,&(2.6b)}$$
where $\epsilon_\si(\mu)\in\{\pm 1\}$. Similar equations
hold for any other affine algebra $g$, and more generally for any RCFT.
The field $K_{r,k}$ here is generated over $\Q$ by all elements
$S_{\la,\mu}$; if instead we are only interested in the permutation $\mu\mapsto
\si\mu$, and not the `parities' $\epsilon_\si(\mu)$, then we are more concerned
with the effective Galois group Gal$(L_{r,k}/\Q)$ coming from the subfield
$L_{r,k}$ generated over $\Q$ by the fusion eigenvalues $\chi_\la(\mu)$. 

Incidentally, Galois orbits tend to be nicely behaved -- see e.g.\ Thm.\ 8 below.
They also have been studied in the `elements of finite order' Lie group
context -- see e.g.\ [18,20].

Galois group considerations are central to many arguments in this paper,
so next we will quickly review the cyclotomic Galois group.
The cyclotomic field $\Q_n:=\Q[\exp[2\pi\i/n]]$ consists of all polynomials
in $\xi_n:=\exp[2\pi\i/n]$. The Galois group ${\cal G}_n:=\,$Gal$(\Q_n/\Q)$ is
isomorphic
to the multiplicative group $(\Z/n\Z)^{\times}$ of integers coprime to $n$,
taken mod $n$. More precisely, any automorphism $\si\in{\cal G}_n$ corresponds
to some integer $\ell\in(\Z/n\Z)^{\times}$, in such a way that $\si\xi_n=\xi_n^\ell$.
We write $\si_\ell$ for this $\si$. The classic example of a Galois
automorphism is complex conjugation, which always corresponds to
$\ell=-1$. A subfield $F$ of $\Q_n$ will have Galois group Gal$(F/\Q)$
isomorphic to a
factor group (equivalently here, a subgroup) of $(\Z/n\Z)^\times$.

The previous properties of $S$ are all well known. The following one, which
relates $S$ entries at fixed points to $S$ entries at both smaller rank and
level, appears to be new. We will call it {\it fixed-point factorisation}.

Let $\varphi$ be a fixed point of $J^d$ for $A_{r,k}$. Then we will show that
$\chi_\la(\varphi)=0$ unless

\item{($*$)} for each $i=1,\ldots,\r/d$, there are precisely
$d$ indices $\ell^{(i)}_1<\cdots<\ell^{(i)}_d$ for which $\la(\ell_j^{(i)})
\equiv -i$ (mod $\r/d$).

\noindent Assume this for now. ($*$) implies ${\r\over d}$ will
divide $t(\la)$ -- which we already know -- but it is much stronger.
Write ${\cal NZ}(d)$ for the set of all weights $\la\in P_+^{r,k}$ which
obey $(*)$. We will see below that
$$\la\in{\cal NZ}(d)\qquad\Longleftrightarrow\qquad\chi_\la\bigl({kd\over \r}
\sum_{i=0}^{\r/d-1}w^{di}\bigr)\ne 0\ .\eqno(2.7a)$$
The fixed-point argument of this last equation has truncated weight $0^\prime$. 

Consider any $\la\in{\cal NZ}(d)$. Let $\pi$ be the unique
permutation of $\{1,\ldots,\r\}$ defined by the following rule:

\itemitem{} \qquad for each $1\le i\le {\r
\over d}$ and $1\le j\le d$, put $\pi(i+(j-1){\r\over d})=\ell^{(i)}_j$.

\noindent $\pi$ will exist iff $(*)$ holds. For each such
$i$, let $\la'{}^{(i)}$ denote the weight in $P_+^{d-1,kd/\r}$ with
Dynkin labels
$$\la'{}^{(i)}_j={\la(\ell^{(i)}_j)-\la(\ell^{(i)}_{j+1})\over \r/d}-1
\eqno(2.7b)$$
for $j=1,\ldots,d-1$. 
As above, let $\varphi'\in P_+^{d-1,kd/\r}$ be the truncated weight $(\varphi_0,
\varphi_1,\ldots,
\varphi_{d-1})$. Then we obtain the `factorisations'
$$\eqalignno{
S_{\la,\varphi}=&\,{\rm sgn}\,\pi\,\,\xi\,(-1)^{t(\la)(1-d/\r)}\bigl({\r\over
\k}\bigr)^{{\r/d-1\over 2}}S'_{\la'{}^{(1)},\varphi'}\cdots
S'_{\la'{}^{(\r/d)},\varphi'}&(2.8a)\cr
\chi_{\la}(\varphi)=&\,{\rm sgn}\,\pi\,\,\xi\,(-1)^{t(\la)(1-d/\r)}
\chi'_{\la'{}^{(1)}}(\varphi')\cdots
\chi'_{\la'{}^{(\r/d)}}(\varphi')&(2.8b)}$$
where $\xi$ is the $\k d/\r$-th root of unity equal to $\exp[2\pi\i\,t'(\varphi'
+\rho')\sum_i(\la_d^{(i)}+i-\r/d)]$, and where
primes denote quantities in $A_{d-1,kd/\r}$ (we take $S'=\chi'=1$
for $d=1$). 

Perhaps some examples at low rank and level will be helpful. For
$r=3,k=4$, the only fixed points are
$(\varphi_0,\varphi_1,\varphi_2,\varphi_3)=(2,0,2,0)$, $(0,2,0,2)$, and
$(1,1,1,1)$. $\NZ(1)$ 
consists of the $J$-orbits of (4,0,0,0) and (2,1,0,1), for a
total of 8 weights out of the full 35. $\NZ(2)$
contains $\NZ(1)$ plus the $J$-orbits of $(3,0,1,0),(2,2,0,0)$ and
$(2,0,2,0)$, increasing the number of weights to 18 out of 35. All 
three fixed points are in the simple-current orbits of weights of
the special type indicated in (2.7a) (for $d=1$ or 2). Therefore, for these
fixed points $\varphi$, we must have $S_{\varphi,\la}\not=0$ for all 
weights $\la$ in the appropriate $\NZ(d)$.   

For $r=3,k=8,d=2$, however, there are fixed points such as $\varphi=(3,1,3,1)$
that are not of the type in (2.7a), i.e.\ $(4,0,4,0)$. In this case, we
find that $S_{(3,1,3,1),\la}\not=0$ for only 48 weights
$\la$, while $\NZ(2)$ has cardinality 75 (and $\|
P_+^{3,8}\|=165$). The large discrepancy here between `48' and `75' is not
surprising and is explained by (2.8): $\chi'_{\varphi'}$ will vanish at a fifth
of the points of $P_+^{2,4}$. Incidentally, the total number of weights
satisfying 
$t(\la)\equiv 0\ ({\rm mod}\ \r/d)$ is 85. This means there are 10 weights
that satisfy the $\r$-ality test necessary for
$\chi_\la(\varphi)\not=0$, yet still have $\chi_\la(\varphi)=0$ for all
$\varphi$.

Condition (*) will become more severe as $r$ and $k$ increase. 
For example, with $r=3$ and $d=2$, the numbers of weights in $\NZ(2)$ compared
with those with even $\r$-ality, compared with those in $P_+^{3,k}$ are:
196, 231, and 455 for $k=12$; and 405, 489, and 969 for $k=16$.

As an example of how `factor weights' $\{\la^{\prime (i)}\}$ are
found, consider the weight $\la=(0,0,1,0,0,1,1,1,0,1,1,0)$ at
$r=11,k=6$. Fix $d=4$. The corresponding partition labels $\{\la(\ell)\}$ are 
$\{17,16,14,13,12,10,8,6,5,3,1,0\}$. Those congruent to $-1$ (mod
$\r/d=3$) are $\{17,14,8,5\}$. From these we find $\la^{\prime(1)}=(1,0,1,0)$,
where the zeroth Dynkin label is set so that the factor weight is at
level $kd/\r=2$. We find $\la^{\prime(2)}=(0,0,0,2)$ and
$\la^{\prime(3)}=(1,1,0,0)$ in similar fashion.

For a more general example, consider any hook $\la=aw^1+w^b$. It will lie in
${\cal NZ}(d)$ iff $\r/d$ divides $a+b$, in which case we find
$$\chi_{aw^1+w^b}(\varphi)=\xi\,(-1)^{a+b+c+(a''+1)(c+a'+1)}
\chi_{(a'-1)w'{}^1+w'{}^{c-a'+1}}(\varphi')\eqno(2.8c)$$
where $c=(a+b)d/\r$ and $a={\r\over d}a'-a''$, for
$a''\in\{1,\ldots,{\r\over d}\}$, and where $\xi=1$ unless $b>\r-\r/d$, in which
case $\xi=\exp[2\pi\i\,t'(\varphi'+\rho')\r/d\k]$. The permutation $\pi$ here
is the product of $c-a'+1$ disjoint $a''$-cycles. In this example, each
$\la'{}^{(i)}=0'$ except for $\la^{(a'')}=(a'-1)w'{}^1+w'{}^{c-a'+1}$.
Equation (2.8c) says that hooks in $P_+^{r,k}$ act like hooks
in $P_+^{d-1,kd/\r}$, when their fusion eigenvalues are restricted to fixed
points of $J^d$. The most interesting special case of (2.8c) is
$$\chi_{w^\ell}(\varphi)=\left\{\matrix{\chi'_{w'{}^{\ell d/\r}}
(\varphi')&{\rm if}\ \r/d\ {\rm divides}\ \ell\cr 0&{\rm otherwise}\cr}
\right.\ .\eqno(2.8d)$$

\medskip\noindent{{\smcap Lemma 3}} (fixed-point factorisation): \quad
{\it Choose any $A_{r,k}$, any
divisor $d$ of gcd$\{\r,k\}$, and any $\la\in P_+^{r,k}$. Then exactly one
of the following holds:

\item{(i)} $S_{\la,\varphi}=\chi_{\la}(\varphi)=0$ for all fixed points
$\varphi$ of $J^d$; or

\item{(ii)} $\la\in{\cal NZ}(d)$ and so $\la$ obeys (2.8a), (2.8b)
for every fixed point $\varphi$ of $J^d$.}\medskip

The leading signs in (2.8) are independent of $\varphi$ and so for our
purposes are of no significance. The phase $\xi$ depends only on $\varphi$ and
will often equal 1. Of course the right side of (2.8b) can be
`linearised' by expanding it out using fusion coefficients. Conversely, it
leads to the curious observation that the fusion coefficients of $A_{r,k}$
can be seen in the fusion eigenvalues of $A_{2r+1,2k}$ evaluated
at fixed points.

At present we do not have formulas of equal generality for the other
affine algebras with simple currents. One would expect that $E_6^{(1)}$
would be related in this way to $G_2^{(1)}$, $E_7^{(1)}$ to $F_4^{(1)}$,
$D_r^{(1)}$ for its vector simple current (i.e.\ the one interchanging $w^0$
and $w^1$, and $w^{r-1}$ and $w^r$) to $C_{r-2}^{(1)}$, etc. Perhaps
an algebraic understanding of these equations can be obtained from the
ideas in e.g.\ [7]. 

To prove equations (2.8), first note that
$$P_\ell[\varphi]=\left\{\matrix{{\r\over d}P'_{\ell d/\r}[\varphi']
&{\rm if}\ \r/d\ {\rm divides}\ \ell\cr 0&{\rm otherwise}\cr}
\right.\ ,\eqno(2.9a)$$
from which we immediately obtain that $H_\ell[\varphi]$ equals
$H'_{\ell d/\r}[\varphi']$
or 0 (for $\r|d\ell$, $\r\not|d\ell$, respectively),
for the `complete' symmetric polynomials $H_\ell:=S_{\ell w^1}$, since (2.4b)
for $\la=\ell w^1$ takes a simple form [8]. We have the determinantal
formula [8]
$$S_{\la}={\rm det}(H_{\la(i)-\r+j})_{1\le i,j\le\r}=\sum_{\sigma}{\rm sgn}\,
\sigma\,\,H_{\la(\sigma 1)-\r+1}\cdots H_{\la(\sigma \r)-\r+\r}\ .\eqno(2.9b)$$
In this formula, $H_0$ identically equals 1, and for negative $\ell$,
$H_\ell$ is identically 0.
Evaluated at the fixed point $\varphi$, this will be a sparse matrix: each
row will have at most $d$ nonzero elements, spaced $\r/d$ entries apart.
If $S_\la[\varphi]\ne 0$, then some product $H_{\la(\sigma 1)-\r+1}[\varphi]
\cdots H_{\la(\sigma\r)-\r+\r}[\varphi]\ne 0$, and thus $\{\ell^{(i)}_1,
\ldots,\ell^{(i)}_d\}=\{\sigma i,\sigma(i+{\r\over d}),\ldots,\sigma(i+\r-{\r\over
d})\}$ for each $i$. This shows that $(*)$ is satisfied, and that the
permutation $\pi$ exists. The sum in (2.9b)
can be restricted to those $\si$ in the coset $(S_d)^{\times(\r/d)}\pi\subset
S_\r$, where the $i$-th factor $S_d$ permutes the indices congruent to $i$
(mod $\r/d$). So (2.9b) can 
now be written as the product of determinants, the $i$-th one of which
corresponds to the weight $\la'{}^{(i)}\in P_+^{d-1,kd/\r}$ (note that (2.4f)
is implicit in (2.7b)), which gives us (2.8b). Equation (2.7a) follows from
(2.8b) and the fact that $({kd\over\r}\sum_i w^{di})'=0'$.
Using the product formula (= Weyl denominator formula) for $S_{0,\mu}$, we
can show
$$S_{0,\varphi}=\left({\r\over \k}\right)^{{\r/d-1\over 2}}(S'_{0',\varphi'})^{
\r/d}\ .\eqno(2.9c)$$
Together with (2.8b), this immediately gives us (2.8a). 

\bigskip\noindent{{\bf 3. Fusion-rank of $A_{r,k}$}}\medskip

The original polynomial realisation [13,15] uses the Cartan fusion-generator
$\Gamma=\{w^1,\ldots,w^r\}$, which works by Lemma 2. 
We can do better.  
 From (2.2a) and Lemma 2, we see that $\R_{r,k}\le {\r\over 2}$, with
$\Gamma=\{w^1,\ldots,w^{\lfloor\r/2\rfloor}\}$, where $\lfloor x\rfloor$ is the
largest integer not larger than $x$. For example, the fusion-rank of $A_{1,k}$
and $A_{2,k}$ equals 1 for all $k$, with $\{w^1\}$ a fusion-generator. This
result for $A_2$ was first obtained in [6], though by a more complicated
argument. We also obtain, from Thm.\ 2(3) below (rank-level duality), the bound
$\R_{r,k}\le {k\over 2}+1$.

We begin by collecting a few simple consequences of the previous comments.
Parts (1) and (3) of Thm.\ 1 are technical facts we will use repeatedly in the
rest of the paper. Thm.\ 1(2) gives a fairly strong lower bound on $\R_{r,k}$.
We give some consequences of Thm.\ 1(4) in the paragraph before Conjecture 1.

\medskip
\noindent{{\smcap Theorem 1}}\ (simple-current constraints):\quad {{\bf (1)}}
{\it Let $\Gamma$ be a fusion-generator, and choose any $\mu\in P_+^{r,k}$.
Let $\Gamma_\mu$ be the set of all $\ga\in\Gamma$ for which $\chi_{\ga}(\mu)
\ne 0$. Let $d={\rm gcd}\{\r,k,t(\ga)|_{\ga\in\Gamma_\mu}\}$ (put $d=\r$ if
$\Gamma_\mu=\emptyset$). Then $\mu$ is a $J^{\r/d}$-fixed point.}

\item{{\bf (2)}} (our best lower bound) {\it Let $\Gamma$ be any fusion-generator. Write out the
prime decomposition $D:={\rm gcd}\{\r,k\}=\prod p_i^{a_i}$, where each prime
$p_i$ is distinct. Then
$$\R_{r,k}\ge\sum a_i\ .$$
If $D\ne \r$, we get the stronger bound 
$$\R_{r,k}\ge 1+\sum a_i\ .$$
More precisely, for each $p_i$, and each $\ell$, $1\le \ell\le a_i$,
there must be some $\ga\in\G\cap {\cal NZ}(\r p_i^\ell/D)$ (see Lemma 3)
with gcd$\{D,t(\ga)\}=D/p_i^\ell$. When $D\ne \r$, there must also be some
$\ga\in\G\cap{\cal NZ}(\r/D)$ whose $\r$-ality $t(\ga)$ is a multiple of $D$.}

\item{{\bf (3)}} {\it Suppose $J^{\r/d}\mu=\mu$ and $J^{\r/d}\nu=\nu$ for some
divisor $d$ of $\r$. Then for any weight $\la$,
$\chi_\la(\mu)=\chi_\la(\nu)\ne 0$
implies $t(\la)\,t(\mu)\equiv t(\la)\,t(\nu)$ (mod $d\,\r$).}

\item{{\bf (4)}} {\it When $\k_1$ is some multiple of $\k_2$, then
$\R_{r,k_1}\ge\R_{r,k_2}$.}

\medskip\noindent{{\it Proof}}\quad (1) Let $\mu$ be a $J^c$-fixed point. Then
from the previous remarks, $c$
 must divide $\r$, and $\r/c$ must divide both $k$ and $t(\ga)$ for each
$\ga\in\G_\mu$. Therefore $c$
must be a multiple of $\r/d$. Moreover $\chi_\ga(J^{\r/d}\mu)=\chi_\ga(\mu)$
for all $\ga\in\Gamma_\mu$ (hence all $\ga\in\Gamma$); since $\Gamma$ is
a fusion-generator this means $J^{\r/d}\mu=\mu$, and hence $c=\r/d$.

(2) We know that for every divisor $d$ of $D$, there are $J^{\r/d}$-fixed points
(more than one, unless
$d=D=\r$). Choose such a fixed point $\varphi$, say. Let
$\Gamma_\varphi$ be as in (1) -- necessarily $\G_\varphi\subseteq {\cal NZ}
(\r/d)$. Then, by (1), 
${\rm gcd}\{\r,k,t(\ga)|_{\ga\in\Gamma_\varphi}\}=d$. So we see  there 
must be a subset $\G_d
\subseteq\G$, namely $\G_d=\G_\varphi$, such that gcd$\{D,t(\ga)
|_{\ga\in\G_d}\}=d$. Note that each
$\G_{D/p_i^\ell}$ must contain some weight $\ga$ with 
gcd$\{D,t(\ga)\}=D/p_i^\ell$ (otherwise different $J^{\r/d}$-fixed points 
would not be distinguished by $\Gamma$).
This gives the first bound. If $\r\ne D$, then there will be several
$J^{\r/D}$-fixed points, and in order for $\G$ to distinguish them, $\G_{D}$
must be nonempty. This gives the second bound.

(3) Let $P_\ell$ be the $\ell$th power sum polynomial (2.4c). From
(2.4d) and (2.5a), $P_\ell[\mu]\ne 0$ requires $d$ to divide $\ell$.
Consider the $m=(m_1,m_2,\ldots)$-th term in $Q_\la$ (see (2.4b)); either 
it will vanish at $\mu$, or $d$ will divide each $\ell$ with $m_\ell\ne
0$. Since $P_\ell[\mu]$ lies in the cyclotomic field  
${\Q}[\exp[2\pi \i\,\ell/\k]]$, we find
that $S_\la[\mu]$ lies in the cyclotomic field ${\Q}[\exp[2\pi \i\,d/\k]]$.
Therefore (2.4a) applied to
$\chi_\la(\mu)=\chi_\la(\nu)\ne 0$ gives us the desired conclusion.

(4) First note that we have the containment
${\k_1\over \k_2}\,(P_+^{r,k_2}+\rho) \subset
P_+^{r,k_1}+\rho$. Moreover, for any weight $\gamma$ we have
$\chi_\gamma^{(1)}({\k_1\over\k_2}(\mu+\rho)-\rho)=\chi_\gamma^{(2)}(\mu)$ for all $\mu\in
P_+^{r,k_2}$, where the superscripts indicate that $k_1$ or $k_2$
should be substituted for $k$ in (1.1b). Suppose $\Gamma^{(1)}$ is a 
fusion-generator for 
$A_{r,k_1}$. Then for any $\mu,\nu\in P_+^{r,k_2}$, if we have 
$\chi_\gamma^{(2)}(\mu)=
\chi_\gamma^{(2)}(\nu)$ for all $\gamma\in \Gamma^{(1)}$, then we know
$\mu=\nu$. Now the $\rho$-shifted action of the affine Weyl group at level
$k_2$ will map any weight $\gamma\in \Gamma^{(1)}$ either to some $\ga'\in
P_+^{r,k_2}$ or onto
the `boundary' of $P_+^{r,k}$. In the former case we get $\chi^{(2)}_\ga(\mu)
=\pm \chi^{(2)}_{\ga'}(\mu)$, for some sign independent of $\mu$.
In the latter case $\chi^{(2)}_\ga(\mu)=0$ for any $\mu$, and can be ignored.
Therefore, the set of weights $\ga'$ in $P_+^{r,k_2}$
obtained in this way from those in $\Gamma^{(1)}$ will be a
fusion-generator for $A_{r,k_2}$. \qed
\medskip

Equation (2.3a) suggests that the fusion-generators for $A_{r,k}$ should be
related to those of $A_{k-1,r+1}$. This is indeed so:

\medskip\noindent{{\smcap Theorem 2}}\ (rank-level duality):\quad 
{{\bf (1)}} {\it Suppose $\r$ does not divide $k$.
Then $\R_{r,k}\ge\R_{k-1,r+1}$. Moreover,
if $\Gamma$\ =\ $\{\ga^1,\ldots,\ga^n\}$ is a fusion-generator for $A_{r,k}$,
then $\widetilde{\Gamma}$\ =\
$\{\tilde{J}^{a_1}\,\tau\ga^1,\ldots,\tilde{J}^{a_n}\, \tau\ga^n\}$ is one for
$A_{k-1,r+1}$, where each $a_i$ is chosen so that gcd$\{a_i\r+t(\ga^i)
,k\}\ =\gcd\{t(\ga^i),\r,k\}$ for each $i$.}

\item{{\bf (2)}} {\it If $\r$ does not divide $k$, and $k$ does not divide $\r$,
then $\R_{r,k}=\R_{k-1,r+1}$;
in this case if $\Gamma$ is a fusion-basis for $A_{r,k}$,
then $\widetilde{\Gamma}$, defined in (1), will be one for $A_{k-1,r+1}$.}

\item{{\bf (3)}} {\it If $\r$ does divide $k$, then $\R_{r,k}\le \R_{k-1,r+1}\le
\R_{r,k}+1$. Using the notation of (1), $\{\tilde{J}\tilde{0},\tau\ga^1,\ldots,
\tau\ga^n\}$ is a fusion-generator for $A_{k-1,r+1}$.}

\medskip\noindent{{\it Proof}}\quad (1) Any weight of $P_+^{k-1,r+1}$
can be expressed as $\tilde J^b\tau\mu$ for some integer $b$ and some weight
$\mu\in  P_+^{r,k}$. So, it suffices to consider any
  $\mu,\mu'\in P_+^{r,k}$ and $b\in\Z$ for which
$$\widetilde{\chi}_{\tilde{J}^{a_i}\,\tau\ga^i}(\tau\mu)= \widetilde{\chi}_{
\tilde{J}^{a_i}\,\tau\ga^i}(\tilde{J}^b\,\tau\mu')\qquad \forall i\ ,\eqno(3.1a)$$
and show that this implies $\tau\mu=\tilde{J}^b\,\tau\mu'$. (3.1a) becomes 
$$\chi_{\ga^i}(\mu)=\exp[2\pi \i\,\{\r a_i+t(\ga^i)\}\,\{t(\mu)-t(\mu')
-\r b\}/\r k]\,\chi_{\ga^i}(\mu')\ .\eqno(3.1b)$$
Define $\G_\mu$ as in Thm.\ 1(1). Because $\r$ does not divide $k$, we know
$\G_\mu\ne \emptyset$. Equation (3.1b) and Thm.\ 1(1) imply that $\mu$ and $\mu'$ will
both be 
$J^{\r/d}$-fixed points, where $d={\rm gcd}_{\ga^i\in\G_\mu}\{a_i\r+t(\ga^i),k
\}$. Then $\tau\mu$ and
$\tilde{J}^b\tau\mu'$ will both be $\tilde{J}^{k/d}$-fixed points.
For each $\ga^i\in\Gamma_\mu$, Thm.\
1(3) and (3.1a) imply
$$\{\r a_i+t(\ga^i)\}\,\{t(\mu)-t(\mu')-\r b\}\equiv 0\quad ({\rm mod}\ d\,k)\
.\eqno(3.1c)$$
For each prime $p|k$, write $p^a$ and $p^{a'}$ for the exact powers dividing
$k$ and $d$, respectively: i.e.\ $p^a\|k$ and $p^{a'}\|d$. So $a\ge a'$. If
$a=a'$, then $p^a$ must divide both $\r$ and $t(\mu)-t(\mu')$, by (2.5b).
If $a>a'$, then $p^{a'}\|(\r a_i+t(\ga^i))$, for some $\ga^i\in\G_\mu$.
Therefore (3.1c) tells us
that $L:=\{t(\mu)-t(\mu')-\r b\}/k$ is an integer.
Equation (3.1b) then implies $\chi_{\ga^i}(\mu)=\chi_{\ga^i}(J^L\mu')$ for all $i$.
Therefore $\mu=J^L\mu'$, so we may take $\mu=\mu'$ in (3.1a), and absorb the
$L$ into $b$. Then $\r/d$
must divide $L$, i.e.\ $k/d$ must divide $b$, i.e.\ $\tilde{J}^b\tau\mu
=\tau\mu$, and we see that (3.1a) can only be trivially satisfied. Hence
$\R_{k-1,r+1}\le \R_{r,k}$.

(2) is immediate from part (1).

(3) The first inequality comes from (1). That the given set is a
fusion-generator follows by the proof of (1). More precisely, by
replacing $\tilde{J}^{a_i}\,\tau\ga^i$ with $\tilde{J}\tilde{0}$ in (3.1a) 
implies $L\in\Z$. The rest of the argument is as before.
\qed\medskip

The Chinese Remainder Theorem tells us that it is always possible to choose
the $a_i$'s in Thm.\ 2(1).
Incidentally, in all cases of which we know,
$\R_{k-1,r+1}=1+\R_{r,k}$ when $\r<k$ divides $k$.
 
Earlier we suggested the upper bound $\R_{r,k}\le \r/2$, and now we also
know $\R_{r,k}\le {k\over 2}+1$ (or $k/2$ if $k$ fails to divide $\r$).
In fact we can do much better than this for most pairs $(r,k)$. The argument
relies on the cyclotomic Galois group ${\cal G}_n$ described briefly in the
previous section.

\medskip\noindent{{\smcap Theorem 3}}\ (Galois considerations):\quad {\it  
$\Gamma_{\div}:=\{w^d\,:\,2d\le \r\ {\rm and}\ d\ {\rm divides}\ \k\,\}$ is a
fusion-generator for $A_{r,k}$, called the {\it divisor generator}. A related
fusion-generator is $\G_{\div}^\tau$, defined by
$$\G_{\div}^\tau:=\left\{\matrix{\{w^d\,:\,2d\le k\ {\rm and}\ d\
{\rm divides}\ \k\}&{\rm when}\ k\ {\rm does\ not\ divide}\ \r\cr
\{w^k\}\cup\{w^d\,:\,2d\le k\ {\rm and}\ d\
{\rm divides}\ \k\}&{\rm when}\ k\ {\rm divides}\ \r\cr}\right.$$
Moreover, each $w^d$ in $\Gamma_{\div}$ and $\Gamma_{\div}^\tau$ can be
replaced with any hook $\ell w^1+ w^{d-\ell}$ for $1\le \ell\le d$.}

\medskip\noindent{{\it Proof}}\quad The key observation here is that, because
each $x_j$ is a $\k$-th root of unity, for any $\ell$ there will exist a Galois
automorphism $\sigma\in {\cal G}_{\k}$ for which
$$\sigma P_d(x_1,\ldots,x_{\r})=P_\ell(x_1,\ldots,x_{\r})\ ,\eqno(3.2)$$
where $d={\rm gcd}\{\ell,\k\}$. 

Suppose, for all $d\le \r/2$ dividing $\k$, that
$$\chi_{w^d}(\mu)=\chi_{w^d}(\mu')  .\eqno(3.3a)$$
We will show this implies $\mu=\mu^\prime$. (3.3a) and (2.4a) give 
$$S_{w^d}[\mu]=\xi^d\,S_{w^d}[\mu']\eqno(3.3b)$$
for all $d\le{\r/ 2}$ dividing $\k$,
where $\xi=\exp[-2\pi \i\,(t(\mu)-t(\mu'))/\r\k]$. 

Equation (2.4b) reads
$$S_{w^\ell}[\mu]={(- 1)^{\ell+1}\over\ell}P_{\ell}[\mu]+{\dot Q}_\ell
(P_1[\mu],\ldots,P_{\ell-1}[\mu])\eqno(3.4a)$$
for some polynomial $\dot Q_{\ell}$ homogeneous in the same sense as $Q_\la$
(and so has no constant term). Let $d$ be the smallest $\ell$ with $P_\ell[\mu]
\ne 0$. Then (3.4a) implies $S_{w^\ell}[\mu]=0$ for all $\ell<d$ and
$S_{w^d}[\mu]=\pm{1\over d} P_d[\mu]\ne 0$, so either $d=\r$ (in which case
$\mu=({k\over\r},
{k\over \r},\ldots,{k\over \r})=\mu'$), or $d\le\r/2$ by (2.2a). But (3.2)
requires $d$ to divide $\k$, if it is to be minimal. Thus (3.3b) holds.
However both $S_{w^d}[\mu]$ and $S_{w^d}[\mu']$ lie in $\Q_{\k/d}$, so $\xi$
must be a $\k$-th root of unity.

We next want to show, by induction on $\ell$, that
$$P_\ell[\mu]=\xi^\ell\,P_\ell[\mu']\eqno(3.4b)$$
for all $\ell\le \r/2$. If we could show this, we would be done, because by
(3.4a) it
would force $\chi_{w^\ell}(\mu)=\chi_{w^\ell}(\mu')$ for all $\ell\le\r/2$,
i.e.\ $\mu=\mu'$.
(3.4b) is clearly true for $P_1=S_1$, using (3.3b) with $d=1$.
By (3.2), it is then true for all $\ell$ with gcd$\{\ell,\k\}=1$. Take any
divisor $d\le \r/2$ of $\k$, and suppose (3.4b) is true for all $\ell<d$.
Using (3.3b), equation (3.4a) means that (3.4b) is true for
$\ell=d$, and hence all $\ell$ with gcd$\{\ell,\k\}=d$. Therefore
(3.4b) is indeed true for all $\ell\le \r/2$, and $\mu=\mu'$.

The above remarks continue to hold if we replace each $w^d$ with any
hook $\ell w^1+w^{d-\ell}$ (of all the weights $\la$ with $t(\la)=d$,
only the hooks have the variable $y_d$ appearing nontrivially in the
corresponding polynomial
$Q_\la(y_i)$ -- see e.g.\ p.51 of [8]). Thm.\ 2 applied to $\G_{\div}$ gives us
the fusion-generator
$\G_{\div}^\tau$ (the hooks $dw^1$ and $J0=kw^1$ can be replaced here with
$w^d$ and $w^k$, respectively).\qed\medskip

In many special cases, most notably Cor.\ 1 and Cor.\ 2 below, we can prove that the
divisor generator $\G_{\div}$ is actually a fusion-{\it basis}. Another example:
suppose gcd$\{\r,k\}=p^\ell$ for some prime $p$, so $\k$ will equal $p^mq$ for
some $m\ge \ell$ and some number $q$ coprime to $p$. If all prime divisors of
$q$ are larger than
$\r/2$, then $\G_{\div}$ will be a fusion-basis, and $\R_{r,k}=\ell+1$ (if
$\r\ne p^\ell$) or $\R_{r,k}=\ell$ (if $\r=p^\ell$). The reason is that here
the lower bound for $\R_{r,k}$ from Thm.\ 1(2) agrees with the upper bound
from Thm.\ 3. A special case of this occurs when both $\r$ and $\k$ are
powers of $p$.

In fact we know of only a few examples (for $\r\le k$) where the divisor
generator is not
a fusion-basis. For $r=4$, for example, we find by computer that the
fusion-rank is 
one for $k=5,9,17$ and 21. 
On the other hand, the computer program tells us that the fusion-rank is
2 for $r=4$ and $k=7,11,13$ and 15. This implies, by Thm.\ 1(4), that whenever
$\k$ is a multiple of 12,16,18 or 20, $\R_{4,k}=2$ and $\G_{\div}$ will be
a fusion-basis.

\medskip\noindent{{\smcap Conjecture 1}}:\quad  {\it At fixed rank $r$, the
divisor generator $\G_{\div}$ 
is a fusion-basis for all sufficiently high levels $k$.} \medskip

For reasons of simplicity, the case of greatest interest is when
$\Gamma=\{w^1\}$ is a fusion-generator. The complete solution to this is a
consequence of this theorem:

\medskip\noindent{{\smcap Theorem 4}}:\quad 
{\it $\Gamma=\{w^1,w^2,\ldots,w^m\}$ is a fusion-generator of $A_{r,k}$ iff}
$\Gamma_{\div}\subseteq\Gamma$ {\it or} $\Gamma_{\div}^\tau\subseteq \Gamma$.

\medskip\noindent{{\it Proof}}\quad `$\Leftarrow$' is immediate from Thm.\ 3. 

`$\Rightarrow$' \quad Suppose we could find a polynomial
$$p(x)=x^{m_1}+\cdots+x^{m_\ell}-x^{n_1}-\cdots-x^{n_\ell}\ ,$$
not identically 0, such that:\smallskip

(a) $\ell< \r$,

(b) $1\le m_1<\cdots<m_\ell< \k$ and $1\le n_1<\cdots<n_\ell< \k$,

(c) $x=\exp[2\pi \i a/\k]$ is a root of $p(x)$, for each $a=1,2,\ldots,m$, and

(d) $\sum_{i=1}^\ell m_i=\sum_{i=1}^\ell n_i$.\smallskip

\noindent Then there would exist weights $\la\ne\mu$ in $P_+^{r,k}$ obeying
$\chi_{w^a}(\la)=\chi_{w^a}(\mu)$ for each $a=1,\ldots,m$ -- in other words,
$\G$ could not in this case be a
fusion-generator. To see this, choose any $\r-\ell$ distinct integers $h_i$
such that $h_1=0$, the remaining $h_i$ obey
 $1\le h_i< \k$, and $\{h_i\}\cap\{m_i\}=\{h_i\}\cap\{n_i\}=\emptyset$. The
$h_i$ and $m_j$ together equal the $\r$ values of $\la(i)$, and the $h_i$
and $n_j$ together equal the $\r$ values of $\mu(i)$. Since $p(x)\not\equiv 0$,
we know $\mu\ne\la$. Condition (c) says that $P_a[\la]=P_a[\mu]$ for all
$a\le m$, and (d) is just the statement that $t(
\la+\rho)=t(\mu+\rho)$. Hence (c) together with induction on (3.4a) is equivalent
to saying $\chi_{w^a}(\la)
=\chi_{w^a}(\mu)$ for those $a$, and we are done.

It is easy to find this polynomial in many cases.
In particular, let $d$ be the largest divisor of $\k$ with $2d\le{\rm min}\{
\r,k\}$, and assume $d>m$. Take $p(x)$ to be $(x^4-x^3-x^2+x)(x^{\k-n}+
x^{\k-2n}+\cdots+x^n+1)$ where $n=\k/d$. Then (c) and (d) are automatically
satisfied. $\ell=2d$ here, so (a) will be satisfied unless $d=\r/2$. Also,
(b) will be satisfied unless $n\le 4$, which can only happen if $d=\r/2=k/2$.

This argument breaks down only when $d=\r/2$. However, when $\r/2$ divides $\k$,
there will be $J^2$-fixed points, and by Thm.\ 1(2) we would require some
$\ga\in\G$ with $t(\ga)$ a multiple of $\r/2$ if $\G$ is to be a
fusion-generator.

The ony remaining way $\Gamma$ could fail to contain $\Gamma_{\div}\cap
\Gamma_{\div}^\tau$ is if simultaneously $k|\r$, $\r\ne k$, and $m<k$. But
then Thm.\ 1(2) applies, and $\Gamma$ would not be able to distinguish the
$J^{\r/k}$-fixed points. \qed

\medskip\noindent{{\smcap Corollary 1}}\ (the first-fundamental generator):
\quad {\it $\Gamma=\{w^1\}$ is a fusion-generator iff both:}

\itemitem{(i)} {\it each prime divisor $p$ of $\k$ satisfies $2p>{\rm min}\,\{\r,\,
k\}$, and}

\itemitem{(ii)} {\it either $\r$ divides $k$, or gcd$\{\r,\,k\}=1$.}\medskip

Incidently, the proof of Thm.\ 4 also implies that at least one weight $\ga$
in any fusion-generator must have $t(\ga)\ge d$, where $d$ is the largest
divisor of $\k$ with $d\le \r/2$ and $d\le k/2$. If this $\ga$ is not
a hook, then in fact $t(\ga)$ would have to be strictly larger than $d$.

\medskip\noindent{{\smcap Corollary 2:}}\quad {\it Some fusion-bases for
$A_{r,k}$ are:}

\item{{$\bullet$}} {\it $\Gamma_{\div}
=\{w^1\}$  for $r=1$ and $2$, $\forall k\ge 1$;

\item{{$\bullet$}} $\Gamma_{\div}=\{w^1\}$ for $r=3$ when $k$ is odd;
$\Gamma_{\div}=\{w^1,\,w^2\}$ for $r=3$ when $k$ is even;

\item{{$\bullet$}} $\Gamma^\tau_{\div}=\{w^1\}$  for $k=1$, $\forall r\ge 1$;

\item{{$\bullet$}} $\Gamma^\tau_{\div}=\{w^1\}$ for $k=2$ and any even $r$; both
$\Gamma=\{J0,w^1\}$ and $\Gamma^\tau_{\div}=\{w^1,
w^2\}$  for $k=2$ and any odd $r>1$;

\item{{$\bullet$}} $\Gamma^\tau_{\div}=\{w^1\}$ for $k=3$ and any $\r$ coprime to
$3$; both $\Gamma=\{J0,w^1\}$ and $\G_{\div}^\tau=\{
w^1,w^3\}$ for $k=3$ and any multiple $\r>3$ of $3$;

\item{{$\bullet$}} $\Gamma^\tau_{\div}=\{w^1\}$ for $k=4$ when $r$ is even;
$\Gamma_{\div}=\{w^1,\,w^2\}$ for $k=4$ when $r\equiv 1$ (mod $4$), $r>4$;
and both
$\Gamma=\{J0,\,w^1,\,w^2\}$ and $\G_{\div}^\tau=\{w^1,w^2,w^4\}$ for $k=4$
when $r\equiv 3$ (mod $4$), $r>4$.}

\medskip
Cor.\ 2 follows immediately from Thm.\ 1(2) and Thm.\ 3. Some of these
fusion-bases are collected in the Table. Cor.\ 2 tells us the fusion-rank
when either $r\le 3$ or $k\le 4$.  

In addition, other fusion-bases are $\Gamma_{\div}=\{w^1\}$ for $r=4$ when
$k$ is even, for $r=5$ when
$k$ is coprime to 6, and $\Gamma^\tau_\div=\{w^1\}$ for $k=6$ when $\r$ is 
coprime to 6;
$\Gamma_{\div}=\{w^1,w^2\}$ for $r=5$ when $k\equiv 2,4$ (mod 6), and 
$\Gamma^\tau_\div=\{w^1,w^2\}$ for
$k=6$ when $r\equiv 1,3$ (mod 6); and $\Gamma_{\div}=\{w^1,w^3\}$
for $r=5$ when $k\equiv 3$ (mod 6), and $\Gamma^\tau_\div=\{w^1,w^3\}$ for $k=6$
when $r\equiv 2$ (mod 6). The simplest cases we do not yet know the answer for
are: $r=4$ when $k$ is odd ($\R\le 2$); $r=5$ when 6 divides $k$ ($\R=2$ or 3);
$k=5$ when $r$ is even ($\R\le 3$); and $k=6$ when 6 divides $\r$ ($\R=3$ or 4).

\topinsert
\centerline{
{\vbox{\offinterlineskip
\hrule
\halign{& \vrule# & \strut $\,\,$ \hfill#\hfill $\,\,$ 
& \vrule# &   #      
& \vrule# & \strut $\,\,$ \hfill#\hfill $\,\,$
& \vrule# & \strut $\,\,$ \hfill#\hfill $\,\,$
& \vrule# & \strut $\,\,$ \hfill#\hfill $\,\,$
& \vrule# & \strut $\,\,$ \hfill#\hfill $\,\,$
& \vrule# & \strut $\,\,$ \hfill#\hfill $\,\,$
& \vrule# \cr
height6pt & \omit && \hskip1pt  && \omit 
&& \omit && \omit && \omit &&
\omit &\cr 
& $r\ \setminus\  k$ && \omit && 1 && 2 && 3 && 4 && 5 & \cr
height3pt & \omit && \omit && \omit && \omit && \omit && \omit && \omit
&\cr \noalign{\hrule}
height1pt & \omit && \omit && \omit && \omit 
&& \omit && \omit && \omit
&\cr \noalign{\hrule}
height6pt & \omit && \omit && \omit 
&& \omit && \omit 
&& \omit && \omit
&\cr & 1 && \omit
&&$\underline{\overline{{\bf\vert}\{w^1\}{\bf\vert}}}$ $\updownarrow$&&
$\{w^1\}$ &&$\{w^1\}$ $\updownarrow$&&$\{w^1\}$ &&$\{w^1\}$ $\updownarrow$&\cr height3pt &
\omit &&\omit && \omit && \omit  && \omit && \omit && \omit &\cr
\noalign{\hrule} height6pt & \omit && \omit && \omit && \omit 
&& \omit && \omit && \omit
&\cr & 2 && \omit
&&${\bf\vert}\{w^1\}{\bf\vert}$ $\updownarrow$&& $\{w^1\}$ $\updownarrow$
&&$\{w^1\}$ &&$\{w^1\}$ $\updownarrow$&&$\{w^1\}$ $\updownarrow$&\cr height3pt & \omit && \omit &&
\omit && \omit && \omit && \omit && \omit &\cr \noalign{\hrule}
height6pt & \omit && \omit && \omit && \omit 
&& \omit && \omit && \omit
&\cr & 3 && \omit
&&${\bf\vert}\{w^1\}$ $\updownarrow$&& $\overline{\{w^1,w^2\}}{\bf\vert}$ 
&&$\{w^1\}$ $\updownarrow$&&$\{w^1,w^2\}$ &&$\{w^1\}$ $\updownarrow$&\cr height3pt & \omit && \omit &&
\omit && \omit && \omit && \omit && \omit &\cr \noalign{\hrule}
height6pt & \omit && \omit && \omit && \omit 
&& \omit && \omit && \omit
&\cr & 4 && \omit
&&$\{w^1\}$ $\updownarrow$&& $\underline{\{w^1\}}$
$\updownarrow$&&$\{w^1\}$ $\updownarrow$&&$\{w^1\}$ $\updownarrow$&&$\{w^2\}$ &\cr height3pt & \omit && \omit &&
\omit && \omit && \omit && \omit && \omit &\cr \noalign{\hrule}
height6pt & \omit && \omit && \omit && \omit 
&& \omit && \omit && \omit
&\cr & 5 && \omit
&&$\{w^1\}$ $\updownarrow$&& $\{w^1,w^2\}$ &&$\overline{\{w^1,w^3\}}$
&&$\overline{\{w^1,w^2\}}$ &&$\{w^1\}$ $\updownarrow$ &\cr height3pt & \omit && \omit &&
\omit && \omit && \omit && \omit && \omit &\cr \noalign{\hrule}
height6pt & \omit && \omit && \omit && \omit 
&& \omit && \omit && \omit
&\cr & 6 && \omit
&&$\{w^1\}$ $\updownarrow$&& $\{w^1\}$ $\updownarrow$&&$\{w^1\}$
$\updownarrow$&&$\{w^1\}$ $\updownarrow$&&$\{w^1,w^2\}$ &\cr height3pt & \omit && \omit && \omit &&
\omit && \omit && \omit && \omit &\cr \noalign{\hrule}
height6pt & \omit && \omit && \omit && \omit 
&& \omit && \omit && \omit
&\cr & 7 && \omit
&&$\{w^1\}$ $\updownarrow$&& $\{w^1,w^2\}$ &&$\underline{\{w^1\}}$
$\updownarrow$&&$\{w^1,w^2,w^4\}$ &&$\{w^1\}$ $\updownarrow$&\cr height3pt & \omit && \omit && \omit &&
\omit && \omit && \omit && \omit &\cr \noalign{\hrule}
height6pt & \omit && \omit && \omit && \omit 
&& \omit && \omit && \omit
&\cr & 8 && \omit
&&$\{w^1\}$ $\updownarrow$&& $\{w^1\}$ $\updownarrow$&&$\{w^1,w^3\}$
&&$\underline{\{w^1\}}$ $\updownarrow$&&$\{2w^2+w^5\}$ $\updownarrow$&\cr height3pt & \omit && \omit &&
\omit && \omit && \omit && \omit && \omit &\cr \noalign{\hrule}
}}}}

\smallskip\smallskip
\leftskip=1cm
\rightskip=1cm
\noindent
\baselineskip=12pt
{\it Table.} Listed are $A_{r,k}$ fusion-bases 
for low ranks and/or levels.
The symbols $\bf\vert$\ in rows of the Table delimit sequences of 
fusion-bases that repeat
indefinitely as the level $k$ increases. For increasing ranks
$r$, overlines and underlines work similarly in the
columns. `$\updownarrow$' signifies that $N_{w^1}$ is invertible (see Section 4).\bigskip
\leftskip=0cm \rightskip=0cm\baselineskip=15pt \endinsert

Obviously to go further we need a better lower bound. Thm.\
1(2) is the best we have,
but it only exploits the presence of fixed points.

\bigskip\noindent{{\bf 4. The fusion matrix of $w^1$}}\medskip

There are many times when
it is useful to know whether particular $S$ matrix elements are nonzero. This
is the case for example in almost every modular invariant partition function
classification attempt -- e.g.\ see the underlying assumption in [17].
It is especially useful to answer this for the first fundamental weight $w^1$
-- in Thm.\ 5 below we give some consequences.

For later convenience, define the sets
$$\eqalignno{{\cal P}_{r,k}:=&\,\{p\ {\rm prime}\ :\ p\le{\rm min}\{\r,k\}\
{\rm and}\ p\ {\rm divides}\ \k\}&(4.1a)\cr
\Z_{\ge}X:=&\,\{\sum_{x\in X}a_xx\ :\ a_x\in \Z_{\ge 0}\}&(4.1b)}$$
where $X$ in (4.1b) is any set of natural numbers. $\Z_{\ge}X$ is the
set of all possible sums (repetitions allowed) of elements of $X$. For example,
 $\Z_{\ge}\{n\}=\{0,n,2n,\ldots\}$ is the set of all nonnegative multiples of
$n$.

\medskip\noindent{{\smcap Theorem 5}}:\quad {{\bf (1)}} {\it Suppose $S_{w^1,
\mu}=0$. Then $S_{\la,\mu}=0$ unless $t(\la)\in\Z_{\ge}{\cal P}_{r,k}$.
Both $k$ and $\r$ must lie in $\Z_{\ge}{\cal P}_{r,k}$.}

\item{{\bf (2)}} {\it Suppose there is only one prime divisor $p$ of $\k$
not larger than min$\{\r,k\}$. Then $S_{w^1,\mu}=0$ iff $\mu$ is a fixed point.}
\medskip 

\noindent{{\it Proof}}\quad  When $k\ge \r$, part (1) follows by considering the polynomial expression (2.4b) and using
the Galois argument of (3.2): $P_\ell[\mu]\ne 0$ requires $\ell\in\Z_{\ge}
{\cal P}_{r,k}$. Taking $\la=J0$ gives us $k\in\Z_{\ge}{\cal P}_{r,k}$,
and $\la={w^\r}$ (see (2.4f)) gives us $\r\in\Z_{\ge}{\cal P}_{r,k}$.
When $k<\r$, to show that we can restrict to
primes $p\le k$, we use rank-level duality (2.3a) to get that $\tilde{t}(\tau\la)
\in\Z_{\ge}{\cal P}_{r,k}$ and then $t(\la)\in\Z_{\ge}{\cal P}_{r,k}$
follows from (2.3b)
and the fact that $k\in\Z_{\ge}{\cal P}_{r,k}$. For part (2), use
part (1) and  Thm.\ 1(1) to get that $\mu$ must be
fixed by $J^{\r/p}$. \qed\medskip

Note that the hypothesis of (2) holds whenever $\k$ is a power of a prime.
This special case follows directly from (4.2) below, by using Gauss' Lemma
on factorising integral polynomials, and evaluating
certain factored polynomials at 1. Thm.\ 5(2) however is much more general.

$N_{w^1}$ is invertible iff $S_{w^1,\mu}\ne 0$ for all $\mu\in
P_+^{r,k}$. Equivalently, $N_{w^1}$ is invertible iff
$$\sum_{j=1}^{\r}\exp[2\pi \i\,\mu(j)/\k]\ne 0\qquad \forall \mu\in P_+^{r,k}\
.\eqno(4.2)$$
It is not hard to show that for $k\le 4$ or $\r\le 4$, $N_{w^1}$ is invertible
iff gcd$\{\r,k\}=1$; in fact, for those $r,k$, $\chi_{w^1}(\mu)=0$ only
for fixed points $\mu$. The identical conclusion holds for many other $r$ and
$k$, as we saw in Thm.\ 5(2). But Thms.\ 6(4),(5) below say that these
cases are uncharacteristically well-behaved. For example, when $\r=5$,
if 6 divides $\k\ge 12$, then $N_{w^1}$ will not be invertible, even though there are
no fixed points.

\medskip\noindent{{\smcap Theorem 6}}\ (invertibility):\quad {\bf (1)} {\it
$N_{w^1}$ is
invertible iff $\tilde{N}_{\tilde{w}^1}$ is, where the latter is the fusion
matrix for $A_{k-1,r+1}$.}

\item{{\bf (2)}} {\it If gcd$\{\r,k\}\ne 1$, then $N_{w^1}$ cannot be invertible.}

\item{{\bf (3)}} {\it $N_{w^1}$ is invertible if either $\r\not\in\Z_{\ge}
{\cal P}_{r,k}$ or $k\not\in\Z_{\ge}{\cal P}_{r,k}$.}

\item{{\bf (4)}} {\it Suppose $pq$ divides $\k$, where $p$ and $q$ are distinct
primes for which $\r\in\Z_{\ge}\{p,q\}$ -- i.e.\ there exist nonnegative
 integers $a,b$ such that $ap+bq=\r$.
If $\k\ge pq(\lceil {a\over q}\rceil +\lceil {b\over p}\rceil)$, then
$N_{w^1}$ will not be invertible ($\lceil x\rceil$ here denotes the smallest
integer not smaller than $x$ -- e.g.\ $\lceil 2\rceil=2$, $\lceil 3.1\rceil=4$).}

\item{{\bf (5)}} {\it Suppose $p_1,p_2,\ldots,p_n$ are primes dividing $\k$
for which $\r\in\Z_{\ge}\{p_1,\ldots,p_n\}$ -- i.e.\ there exist nonnegative
 integers $a_i$ such that $\sum a_ip_i=\r$.
If $\k\ge p_ip_j\sum_{h=i}^j a_h$ for any $i<j$, then $N_{w^1}$ will not be
invertible.}

\medskip\noindent{{\it Proof}}\quad (1) follows directly from (2.3a).
(2) exploits the fact (see (2.2b)) that $\chi_{w^1}(\varphi)=0$ for any fixed
point $\varphi$. (3) is a corollary of Thm.\ 5(1).

(4) We want to construct a particular $\mu\in P_+^{r,k}$ such that
$\chi_{w^1}(\mu)=0$. To do this we find an arithmetic sequence ${\k\over
p}\Z+c_i$ for each $i=1, \ldots,a$, and an arithmetic sequence ${\k\over
q}\Z+c_j'$ for each $j=1, \ldots,b$, such that none of these $a+b$ sequences
intersect. This is easy to do, provided $\k$ is big enough. Choose as the
$c_i$'s $0,{\k\over q}, \ldots,{\k\over q}(q-1),1,1+{\k\over
q},\ldots,1+{\k\over q}(q-1)$, etc, until we have chosen $a$ of them (the last
one will be $\lceil {a\over q}\rceil-1$ plus some multiple of ${\k\over q}$).
Next choose as the $c_j'$'s $\lceil{a\over q}\rceil,\lceil{a\over
q}\rceil+{\k\over p},\ldots$, until we have chosen $b$ of them (the last one
will be $\lceil{a\over q}\rceil +\lceil{b\over p}\rceil -1$ plus some multiple
of ${\k\over p}$). Our $a+b$ sequences will be disjoint, provided the bound on
$\k$ is satisfied, and will intersect the interval $0\le x<\k$ in precisely
$ap+bq=\r$ points. Let $\mu$ be the unique weight in $P_+^{r,k}$ whose
$\mu(\ell)$ equal those $\r$ points. Then $\chi_{w^1}(\mu)=0$, because the sum
in (4.2) along each of the $a+b$ sequences is 0.

(5) follows immediately from similar considerations: we are looking for
$a_i$ series ${\k\over p_i}\Z+c_{ij}$, where $c_{ij}\not\equiv c_{i\ell}$
(mod ${\k\over p_i}$) for $j\ne \ell$, and $c_{ij}\not\equiv c_{h\ell}$
(mod ${\k\over p_i p_h}$) for $i\ne h$. The choice $c_{ij}=
j-1+\sum_{\ell=1}^{i-1}a_\ell$ works.\qed \medskip

The proofs of Thms.\ 6(4),(5) are constructive: their zeros arise when (4.2)
finds itself a sum of
terms such as $\sum_{a=1}^p\xi^a$ for $\xi$ a primitive $p$-th
root of unity. A simple example of Thm. 6(4) is at
$\r=11,\k=30$. With $p=3,q=5$, and $a=2,b=1$, the bound is saturated. One
finds $c_1=0,c_2=6$ and $c_1^\prime=1$. These yield $0,10,20$; $6,16,26$; and
$1,7,13,19,25$; respectively. So, there is a zero for the weight given by
$\{\mu(1),\ldots,\mu(\r)\}$ =
$\{26,25,20,19,16,13,10,7,6,1,0\}$.\medskip  

\noindent{{\smcap Conjecture 2}}:\quad {\it For $A_{r,k}$,
$N_{w^1}$ fails to be invertible
iff one can find distinct primes $p_i\le{\rm min}\{\r,k\}$ dividing
$\k$ and nonnegative integers $a_i,b_i$ such that $\r=\sum_{i}a_ip_i$ and
$k=\sum_ib_ip_i$.}\medskip

In other words, we conjecture that the condition of Thm.\ 6(3) is an `iff'.
Note that one way this condition will be satisfied is if gcd$\{\r,k\}\ne 1$.
The conditions in Thms.\ 6(4),(5) are strongest when we take $r<k$ (which
without loss of generality we can). Also, the bound in 6(5) is best when the
$p_i$ are labelled so that the largest are given indices near $n/2$.
In practice the most useful special
case of Thms.\ 6(4),(5) is: If one can find an odd prime $p\le \r$ for which $2p$ divides
$\k$ and $k\ge 3p-1$, then $N_{w^1}$ will not be invertible.
The analogue of Thm.\ 1(4) is also valid here, but is not very useful.

The answer to Question 2 for small $r$ and $k$ is indicated in the table.
Computer checks were performed for $r\le 9$ and all levels $k>r$ such that
$\dim P_+^{r,k}$ $< 300,000$. The results were consistent with 
Conjecture 2. Conjectures 1 and 2 are
the simplest guesses consistent with our results, but it would be nice to
test them against additional numerical data.

Incidentally, conditions like `$\ell\in\Z_{\ge}\{n_1,\ldots,n_m\}$'
are only strong when $\ell$ is small. For example, given any
coprime numbers $m$ and $n$, there are only $(m-1)\,(n-1)/2$ positive
integers $\ell$ which do not lie in $\Z_{\ge}\{m,n\}$ --
the largest such $\ell$ is $mn-m-n$.
So for fixed $\r$, we know {\it Conjecture 2 will hold for all sufficiently
large $k$}.

\bigskip\noindent{{\bf 5. Extensions}}\medskip

Because the fundamental weights are much simpler, the most interesting
fusion-generators are the ones which consist only of fundamental weights:
$\Gamma\subseteq\{w^1,\ldots,w^r\}$. We can speak of {\it
fundamental-fusion-generators} and {\it fundamental-fusion-rank} ${\cal FR}_{
r,k}$. All of the results in Sections 3 and 4 also apply directly to ${\cal
FR}_{r,k}$. By definition, ${\cal R}_{r,k}\le {\cal FR}_{r,k}$, and
Conjecture 1 predicts that, for fixed $r$, ${\cal R}_{r,k}={\cal FR}_{r,k}$
for all sufficiently large $k$. Note however from the Table that ${\cal FR}_{
8,5}={\cal FR}_{4,9}=2$ while ${\cal R}_{8,5}={\cal R}_{4,9}=1$.

Because of (2.8d), we can strengthen here the bound in Thm.\ 1(2). For example,
if ${\cal FR}_{r,k}$ equals the bound given in Thm.\ 1(2), then so must
${\cal FR}_{\r/d-1,k/d}$ for all divisors $d$ of gcd$\{\r,k\}$.

One can also ask Question 2 for other weights, most importantly the
other fundamental weights, and again (2.8b) will be very useful. For
example, we know $\chi_{w^2}$ will vanish at some $J^5$-fixed point of
$A_{9,14}$, because $N_{w^1}$ is not invertible for $A_{4,7}$.

Of course Questions 1 and 2 can and should be asked of the fusion algebras for the
other affine algebras, and similar arguments will apply. We have not
investigated them, except to find some fusion-bases for $C_{2,k}$ and $G_{2,k}$
on the computer, and to get Thm.\ 7 below for $G_{2,k}$.
Of course $\R_k(C_2^{(1)})$ must equal 2 for
any even $k$, and we find the rank is also 2 for all odd $k<26$ (the limit
of our computer check), save $k=1,3$ and 9. For $k=1$ and 9, the only
fusion-bases are
$\{w^1\}$ and $\{2w^1+6w^2\}$, respectively. At $k=3$ there are four
different fusion-bases: $\{2w^1\}$, $\{w^2\}$, $\{2w^1+w^2\}$, and
$\{2w^2\}$. A very tempting conjecture is that the rank $\R(C_{r,k})$ equals 2
for all sufficiently large $k$ (and probably for all $k>9$). The situation for
$G_{2,k}$ however is more surprising:

\medskip\noindent{{\smcap Theorem 7}}:\quad {\bf (1)} {\it When the level $k$
is odd, $\{w^2\}$ is a fusion-basis for $G_{2,k}$.}

\item{{\bf (2)}} {\it $N_{w^2}$ fails to be invertible for $G_{2,k}$ iff either
4 or 30 divides $\k:=k+4$.}\medskip

\noindent{\it Proof}\quad The key here is to reduce the $G_{2,k}$
quantities to $A_{2,k+1}$ quantities, and use the fact that $\{w^1\}$ is a
fusion-basis for $A_{2,k+1}$.

Using (1.1b) and the simple Lie subalgebra $A_2\subset G_2$, we find 
$$\chi_{w^2}(\mu)=\underline{\chi}_{\underline{w}^1}(\underline{\mu})+
	\underline{\chi}_{\underline{w}^2}(\underline{\mu})+1\ ,\eqno(5.1)$$
where underlines denote $A_{2,k+1}$ quantities, and $\underline{\mu}=\mu_1\underline{w}^1+(\mu_1+\mu_2+1)\underline{w}^2$.
So part (1) reduces to the following statement\footnote{$^4$}{{\smal For the
remainder of the proof of part (1), we will switch to $A_{2,k+1}$ notation.}}
 for $A_{2,k+1}$: for any
$\la,\mu\in P_+^{2,k+1}$ with $\la\ne C\la$ and $\mu\ne C\mu$ (only these
nonselfconjugate weights correspond to $G_{2,k}$ ones), does the equality
$$\eqalignno{\cos&(2\pi{\la_1+2\la_2+3\over 3\k})+\cos(2\pi{\la_2-\la_1\over 3\k})+
\cos(2\pi{2\la_1+\la_2+3\over 3\k})&(5.2a)\cr &=
\cos(2\pi{\mu_1+2\mu_2+3\over 3\k})+\cos(2\pi{\mu_2-\mu_1\over 3\k})+
\cos(2\pi{2\mu_1+\mu_2+3\over 3\k})&}$$
force either $\la=\mu$ or $\la={C}\mu$? 

Write $c_1,c_2,c_3$ for the three cosines on the left side of (5.2a), and
write $c_1',c_2',c_3'$ for those on the right. Then (5.2a) says $c_1+c_2+c_3=
c_1'+c_2'+c_3'$, and since $(2\nu_1+\nu_2+3)+(\nu_2-\nu_1)=\nu_1+2\nu_2+3$,
we also get $c_1^2+c_2^2+c_3^2=1+2c_1c_2c_3$ and $c_1'{}^2+c_2'{}^2+c_3'{}^2=1+
2c_1'c_2'c_3'$.

Hit both sides of (5.2a) with the Galois automorphism $\si_2$ (see Section 2).
Since $\cos(2x)=2\cos^2(x)-1$, we obtain
$$c_1^2+c_2^2+c_3^2=c_1'{}^2+c_2'{}^2+c_3'{}^2\ .\eqno(5.2b)$$
Thus any symmetric polynomial in $c_1,c_2,c_3$ will equal the corresponding
symmetric polynomial in $c_1',c_2',c_3'$. In particular 
$$\eqalignno{\bigl(\sin&(2\pi{\la_1+2\la_2+3\over 3\k})
-\sin(2\pi{\la_2-\la_1\over 3\k})-
\sin(2\pi{2\la_1+\la_2+3\over 3\k})\bigr)^2&(5.2c)\cr&=\bigl(
\sin(2\pi{\mu_1+2\mu_2+3\over 3\k})-\sin(2\pi{\mu_2-\mu_1\over 3\k})-
\sin(2\pi{2\mu_1+\mu_2+3\over 3\k})\bigr)^2\ .&}$$
In other words, we know from (5.2a) that the real parts of
${\chi}_{{w}^1}({\la})$ and ${\chi}_{{w}^1}({\mu})$ are equal, and from
(5.2c) that their
imaginary parts are also equal, up to a sign. Hence either ${\la}=
{\mu}$ or ${\la}={C}\mu$, and we have proven part (1).

For part (2), note that $\chi_{w^2}(\mu)=0$ is equivalent to (see (5.1))
$$c_1+c_2+c_3=-{1\over 2}\ ,\eqno(5.3a)$$
 in the above notation. Consider first $k$ odd. Then hitting (5.3a) with
the Galois automorphism $\sigma_2$ gives us $c_1^2+c_2^2+c_3^2
={5\over 4}$, and hence $c_1c_2c_3={1\over 8}$.
We can solve these equations, and we find $8c_i^3+4c_i^2-4c_i-1=0$, i.e.\
$\{c_1,c_2,c_3\}=\{\cos(2\pi{1\over 7}),\cos(2\pi{2\over 7}),\cos(2\pi{3\over
7})\}$. However, these cosines cannot be realised by a weight in $P_+^k(
G_2^{(1)})$.

Next, suppose $k\equiv 2$ (mod 4). We may assume (using $G_{2,k}$ notation)
that exactly two of
the arguments $\{3\mu_1+2\mu_2+5,\mu_2+1,3\mu_1+\mu_2+4\}$ are odd, otherwise
they would all be even and the argument would reduce to the $k$ odd one.
Here we use the automorphism $\sigma_{3\k/2-2}$
and find (relabeling the $c_i$ if necessary) that $c_3^2-c_1^2-c_2^2=-{3\over 4}$.
We can solve for $c_i$ as before, and we find that either $c_3=\cos(2\pi{1\over 5})$
and $\{c_1,c_2\}=\{\cos(2\pi{7\over 30}),\cos(2\pi{13\over 30})\}$, or
$c_3=\cos(2\pi{2\over 5})$
and $\{c_1,c_2\}=\{\cos(2\pi{1\over 30}),\cos(2\pi{11\over 30})\}$.
Either possibility requires 30 to divide $\k$, in order to be realised by a
weight of $G_{2,k}$.
When 30 divides $\k$, we do indeed get  zeros: $\mu=(\k/3-1,\k/30-1,3\k/5-1)\in
P_+^k(G_{2}^{(1)})$ works.

Finally, suppose 4 divides $k$. Then $\mu=(k/4,k/4,k/4)\in P_+^k(G_{2}^{(1)})$ works.
\qed\medskip

(By $w^2$ here we mean the Weyl-dimension 7 fundamental weight of $G_2$,
corresponding to the short simple root.)
However, $\{w^2\}$ will not be a fusion-generator when $k>4$ is even.
Our computer program tells us that for $k\le 24$, the fusion-rank is 1
except for $k=6,12,16$ and 20 (of course this implies it will also be 2
whenever $k+4$ is a multiple of 10, 16, or 24).

\bigskip
\noindent{{\bf 6. Number fields associated with $S$}}\medskip

By the field $K_{r,k}$ we mean the smallest field containing the rationals
and all of the entries $S_{\la,\mu}$ of $S$. Similarly, by the field $L_{r,k}$
we mean
the smallest field containing $\Q$ and all of the values $\chi_\la(\mu)$.
Because of their role in the Galois symmetry (2.6), it is natural to
try to identify these fields. This question was posed in [4], and related
questions have been considered in e.g.\ [18,20]. Another
reason the question is interesting is that, as we shall see, it has a simple
answer! We will
give this answer in Cor.\ 3 below, for the most important case: $A_{r,k}$.

The matrix $S$ for any nontwisted affine algebra $g$ is given in e.g.\ [16].
The expression for
$S_{\la,\mu}$ consists of a sum $s(\la,\mu)$ over the Weyl group of 
$\overline{g}$,
multiplied by a constant $c$. 
For $A_{r,k}$, $s(\la,\mu)$ manifestly lies in the field $\Q_{\r\k}$, 
and 
$$c={\i^{r(r+1)/2}\over \k^{r/2}\,\sqrt{r+1}}\ .$$
Using Gauss sums, which express square-roots of integers as sums of roots
of unity, it can be shown that the constant $c$ lies in either $\Q_{\r}$
if $r$ is even, or $\Q_{\r\k}$ if either $r\equiv 3$ (mod 4) or $k$ is even, or
$\Q_{\r\k}[\sqrt{\pm 2}]$ if both $k$ is odd and $k\r\equiv \pm 2$ (mod 8).
Thus we know $L_{r,k}$ is always a subfield of $\Q_{\r\k}$, and $K_{r,k}$
is always a subfield of $\Q_{4\r\k}$.

Write $[\la]$ for the orbit $\{J^i\la\}$ of $\la$ by the simple currents.
We will find our fields by first computing some Galois orbits. This result
should be of independent value.

\medskip\noindent{{\smcap Theorem 8}}:\quad {\it Consider any $k>2$ and $r
\ne 1$.}

\item{{\bf (1)}} {\it Choose any fundamental weight $w^m$ with $m\le{\rm min}
\{\r-2,k-2\}$, and any Galois
automorphism $\si_\ell$. Then (with one exception)
$\si_\ell w^m\in[w^m]\cup[Cw^m]$ iff $\ell\equiv
\pm 1$ (mod $\k)$; for all other $\ell$ the quantum-dimension $S_{\si_\ell w^m,
0}/S_{0,0}$ of $\si_\ell w^m$ will be strictly greater than that of $w^m$. (The
one exception is $w^2$ for $A_{3,4}$, where each $\sigma_\ell$ fixes $w^2$.)}

\item{{\bf (2)}} {\it When $r\not\equiv 1$ (mod 4), $\si_\ell w^1=w^1$ iff $\ell=1$ (mod $\r\k$).
When $r\equiv 1$ (mod 4) and $k$ is even, then
$\si_\ell w^1=w^1$ iff $\ell=1$ (mod $\r\k/2$).}

\medskip\noindent {\it Proof} \quad (1)
Because of (2.1a), we may assume $m\le \r/2$.
Assume first that $k\ge \r$. From the Weyl denominator formula, we compute
$${S_{\si_\ell w^m,0}\over S_{w^m,0}}=\prod_{n=1}^m{|\sin(\pi \ell n/\k)|^{r-n}
\over\sin(\pi n/\k)^{r-n}}\prod_{n=m+1}^{\r-m}{|\sin(\pi\ell n/\k)|^{\r-n}
\over \sin(\pi n/\k)^{\r-n}}\prod_{n=\r+1-m}^\r{|\sin(\pi\ell n/\k)|^{\r+1-n}
\over\sin(\pi n/\k)^{\r+1-n}}\eqno(6.2a)$$
where we drop the middle product if $m=\r/2$. 
We want to know when (6.2a) equals 1. This is easy, for $k>r\ge 2$, since
$\sin(\pi/\k)<\sin(2\pi/\k)<\cdots<\sin(\pi\r/\k)$. Consider first $m<\r/2$:
of all possible choices
of integers $1\le n_1<n_2<\cdots<n_{r+1}\le \k/2$, the minimum possible
product of $r-1$ $\sin(\pi n_1/\k)$'s, $r-2$ $\sin(\pi n_2/\k)$'s, ...,
$r-m$ $\sin(\pi n_m/\k)$'s, $r-m$ $\sin(\pi n_{m+1}/\k)$'s, ..., $m$ $\sin(\pi
n_{\r-m}/\k)$'s, $m$ $\sin(\pi n_{\r+1-m}/\k)$'s, ..., and 1 $\sin(\pi n_{\r}
/\k)$, is the choice $n_1=1$, $n_2=2$, ..., $\{n_m,n_{m+1}\}=\{m,m+1\}$, ...,
$n_{m+2}=m+2$, ..., $\{n_{\r-m},n_{\r+1-m}\}=\{\r-m,\r+1-m\}$, ..., $n_{m+1}=
m+1$. This immediately forces $\ell\equiv \pm 1$ (mod $\k$) (for $m>1$, just
look at the first term; when $m=1$, $\ell\equiv\pm 2$ is eliminated by seeing
what happens to the second term).

If instead $m=\r/2$, the exponents of $\sin(\pi n/\k)$ in (6.2a)
are no longer nonincreasing: near $n=m+1$ we get the subproduct
$$\cdots\sin(\pi (m-1)/\k)^{\r-m}
\sin(\pi\,m/\k)^{r-m}\sin(\pi\,(m+1)/\k)^{\r-m}\sin(\pi\,(m+2)/\k)^{r-m}
\cdots$$
For $m>2$, the proof that (6.2a) will always be greater than 1 for
$\ell\not\equiv \pm 1$ (mod $\k$), follows from the simple observation that
$\sin(\pi/\k)\,\sin(\pi\,(m+1)/\k)<\sin(2\pi/\k)\,\sin(\pi\,m/\k)$: the
least-harmful place to move `1' to is `2', and the best place to move `$m+1$'
to is `$m$', and yet even that (forgetting the other terms, which will make
matters worse) will increase the product. The remaining case $m=2$ corresponds
to $\r=k=4$, i.e.\ to the given exception.  

This completes the argument for $k\ge \r$. When $k<\r$, apply rank-level
duality (2.3a): it is an exact symmetry of quantum-dimensions, and maps
$J$-orbits to $\tilde{J}$-orbits. $\tau w^m=m\tilde{w}^1$, so we are interested
in the ratio
$${\tilde{S}_{\si_\ell m\tilde{w}^1,0}\over \tilde{S}_{m\tilde{w}^1,0}}=
\prod_{n=1}^{k-2}{|\sin(\pi \ell n/\k)|^{k-1-n}
\over\sin(\pi n/\k)^{k-1-n}}\ \prod_{n=1}^{k-1}{|\sin(\pi\ell \,(n+m)/\k)|
\over \sin(\pi\,( n+m)/\k)}\ .\eqno(6.2b)$$
The rest of the argument is as before: again $m=\r/2$ causes minor problems.

Now consider any $\ell=(-1)^a+b\k$. Applying (2.6b) to the Cartan generators
$\la\in\{w^1,\ldots,w^r\}$ and using (2.4e), we find
$$\si_\ell \mu=C^a\,J^{b\,t(\mu+\rho)}\mu\eqno(6.2b)$$
whenever $\sigma_\ell\in{\rm Gal}(L_{r,k}/\Q)$. Applying (6.2b) to
$\mu=w^1$ gives us part (2). \qed\medskip

\noindent{{\smcap Corollary 3}}:\quad{\it When both $k>2$ and
$r\ne 1$, then $L_{r,k}=\Q_{\r\k}$ and}
$$K_{r,k}=\left\{\matrix{\Q_{\r\k}&{\rm if\ either}\ r\not\equiv 1\ ({\rm
mod}\ 4)\ {\rm or}\ k\ {\rm is\ even}\cr \Q_{\r\k}[\sqrt{\pm 2}]&{\rm if}\
r\ {\rm is\ odd\ and}\ \r k\equiv \pm 2\ ({\rm mod}\ 8)\cr}\right.$$

The proof of the Corollary is immediate from Thm.\ 8, by regarding Galois
orbit sizes: when $r\not\equiv 1$ (mod 4), the Galois orbit of $w^1$ alone
suffices, but when $r\equiv 1$ (mod 4) and $k$ is even, we have $\si_{1+\r\k/2}w^1=w^1$, so
also use $\si_{1+\r\k/2}0=J^{\r/2}0\ne 0$, which is obtained from (6.2b).
What we find in all cases is that for any $\ell\in(\Z/\r\k\Z)^\times$,
$\ell\ne 1$, either $\si_\ell w^1\ne w^1$ or $\si_\ell 0\ne 0$. This tells us
$L_{r,k}=\Q_{\r\k}$, and $K_{r,k}$ is then obtained by adjoining the
constant $c$ shown above.

Similar statements to Thm.\ 8 can be found for other weights. For example,
by rank-level duality the identical result to Thm.\ 8(1) holds for any
$mw^1$, $0\le m\le{\rm min}\{\r-2,k-2\}$, and we can expect similar results for
other hooks.
When $r\equiv 1$ (mod 4) and $k$ is odd, $\Q_{4\r\k}$ is a degree 2 extension
of $K_{r,k}$, which is in turn a degree 2 extension of $\Q_{\r\k}$.

The results
corresponding to Cor.\ 3 for $k=1,2$ or $r=1$ can be easily found, but are more
complicated and hence less interesting. We include them here for
completeness.

\item{$\bullet$} $L_{r,1}=\Q_{\r}$.  $K_{r,1}$ will equal either
$\Q_{\r}$, $\Q_{\r}[\i]$, or $\Q_{\r}[\sqrt{\pm 2}]$, depending on whether
or not $\r\equiv 0,1$ (mod 4), or $\r\equiv 3$ (mod 4), or $\r\equiv\pm 2$
(mod 8), respectively.

\item{$\bullet$} $L_{1,k}=\Q[\cos(2\pi/\k)]$ if $k$ is odd, and $\Q[\cos(\pi/
\k)]$ if $k$ is even. $K_{1,k}$ will equal either $L_{1,k}$, or $L_{1,k}[
\sqrt{2}\sin(2\pi/\k)]$, or $L_{1,k}[\sqrt{2}]$, depending
on whether $k\equiv 0,2$, or $k\equiv 3$, or $k\equiv 1$ (mod 4), respectively.

\item{$\bullet$} $L_{r,2}=\Q_{\r}[\cos(2\pi/\k)]$ if $\r$ is odd, 
 and $\Q_{\r\k}$ if $\r$ is even.
$K_{r,2}$ will equal $L_{r,2}$, unless $\r\equiv 3$ (mod 4) when $K_{r,2}=
\Q_{\r\k}$.

\vfill \eject
\noindent {\bf Acknowledgements}

T.G. thanks A.\ Coste for showing him Questions 1 and 3, and C.\ Cummins
for discussions.
M.W. thanks the High Energy Physics group of DAMTP for hospitality, 
and W. Eholzer for reading the manuscript.  

\bigskip
\noindent{\bf REFERENCES}
\bigskip
\item{1.} Aharony, O.: Generalized fusion potentials, {\it Phys.\ Lett.}
{\bf B306} (1993), 276-282.

\item{2.} Altschuler, D.,  Bauer, M., and Itzykson, C.: The branching
rules of conformal embeddings, {\it Commun.\ Math.\ Phys.}\ {\bf 132} (1990),
349-364.

\item{3.} Bourbaki, N.: {\it Groupes et Alg\`ebres de Lie, Chapitres IV-VI,}
Hermann, Paris, 1968.

\item{4.} Buffenoir, E., Coste, A., Lascoux, J., Buhot, A., and Degiovanni,
P.: Precise study of some number fields and Galois actions occurring in
conformal field theory, {\it Annales de l'I.H.P.: Phys.\ Th\'eor.} {\bf 63}
(1995), 41-79.

\item{5.} Coste, A.\ and Gannon, T.: Remarks on Galois symmetry in
rational conformal field theories, {\it Phys.\ Lett.}\ {\bf B323} (1994),
316-321.

\item{6.} Di Francesco, P.\ and Zuber, J.-B.: Fusion Potentials I,
{\it J.\ Phys.}\ {\bf A26} (1993), 1441-1454.

\item{7.} Fuchs, J., Schellekens, B., and Schweigert, C.: From Dynkin diagram
symmetries to fixed point structures, {\it Commun.\ Math.\ Phys.}\ {\bf
180} (1996), 39-97.

\item{8.} Fulton, W.\ and Harris, J.: {\it Representation Theory: A
First Course}, Springer-Verlag, New York, 1991.

\item{9.} Gannon, T.: Symmetries of the Kac-Peterson modular matrices
of affine algebras, {\it Invent.\ Math.}\ {\bf 122} (1995), 341-357.

\item{10.} Gannon, T.: Kac-Peterson, Perron-Frobenius, and the
classification of conformal field theories, e-print q-alg/9510026,
(1995). 

\item{11.} Gannon, T., Ruelle, Ph., and Walton, M.\ A.:  Automorphism
modular invariants of current algebras, {\it Commun.\ Math.\
Phys.}\ {\bf 179} (1996), 121-156.

\item{12.} Georgieu, G.\ and Mathieu, O.: Cat\'egorie de fusion pour
les groupes de Chevalley, {\it C.\ R.\ Acad.\ Sci.\ Paris} {\bf 315} (1992),
659-662.

\item{13.} Gepner, D.: Fusion rings and geometry, {\it Commun.\ Math.\
Phys.}\ {\bf 141} (1991), 381-411.

\item{14.} Goodman, F.\ and Nakanishi, T.:  Fusion algebras in
integrable systems in two dimensions, {\it Phys.\ Lett.}\ {\bf B262} (1991),
259-264.

\item{15.} Goodman, F.\ and Wenzl, H.: Littlewood-Richardson coefficients
for Hecke algebras at roots of unity, {\it Adv.\ Math.}\ {\bf 82} (1990) 244-265.

\item{16.} Kac, V.\ and Peterson, D.: Infinite-dimensional Lie algebras, theta
functions and modular forms, {\it Adv.\ Math.}\ {\bf 53} (1984), 125-264.

\item{17.} Kreuzer, M.\ and Schellekens, A.\ N.: Simple currents versus
orbifolds with discrete torsion -- a complete classification, {\it Nucl.\
Phys.}\ {\bf B411} (1994), 97-121.

\item{18.} Moody, R.\ V.\  and Patera, J.:  Characters of elements of finite
order in Lie groups, {\it SIAM J.\ Alg.\ Disc.\ Meth.} {\bf 5} (1984), 359-383.

\item{19.} Pasquier, V.\ and Saleur, H.: Common structures between finite 
systems and conformal field theories through quantum groups, {\it Nucl.\ Phys.}
 {\bf B330} (1990), 523-526.

\item{20.} Pianzola, A.: The arithmetic of the representation ring and elements
of finite order in Lie groups, {\it J.\ Algebra} {\bf 108} (1987), 1-33.
 
\item{21.} Verlinde, E.: Fusion rules and modular transformations in 2D
conformal field theory, {\it Nucl.\ Phys.}\ {\bf B300} (1988),
360-376.

\item{22.} Witten, E.: The Verlinde algebra and the cohomology of
the Grassmannian, In: {\it Geometry, Topology and Physics}, Conf.\ Proc.\
and Lecture Notes in Geom.\ Topol.\ vol.\ VI (1995), 357-422.

\end